\documentclass[11pt]{article}
\usepackage{amssymb}
\usepackage{color}
\usepackage{amsfonts}
\usepackage{latexsym}
\usepackage[dvips]{graphics}
\usepackage{epsf, epsfig}
\usepackage{mathdots}

\input prepictex
\input pictex
\input postpictex

\newtheorem{lem}{Lemma}[section]
\newtheorem{theo}{Theorem}[section]
\newtheorem{pro}{Proposition}[section]
\newtheorem{cor}{Corollary}[section]

\newtheorem{prob}{Problem}[section]

\definecolor{orange}{RGB}{255,102,102}
\definecolor{gray}{RGB}{100,100,100}
\definecolor{purple}{RGB}{153, 0, 255}
\definecolor{brown}{RGB}{139, 58, 58}
\definecolor{pink}{RGB}{255,0,204}

\def \cirs {\makebox(0,0)[l]{\circle*{6}}}

\newcommand{\proof}
{{\noindent {\em Proof}.\quad}\setcounter{countclaim}{0}
\setcounter{countcase}{0}}
\newcommand{\proofend}{{\hfill$\Box$}}

\newcommand{\sgap}{\vspace{0.3 cm}}

\newcounter{countfig}

\newcounter{countclaim}

\newcounter{countcase}

\newcommand{\beeq}{\begin{equation}}
\newcommand{\eneq}{\end{equation}}

\newcommand{\beeqn}{\begin{eqnarray*}}
\newcommand{\eneqn}{\end{eqnarray*}}

\def \setb{{\cal B}}

\def \sett{{\cal T}}
\def \N {{\mathbb N}}

\def \GP {{\cal G}{\cal P}}
\def \BP {{\cal BP}}
\def \BPX {{\cal BP}_X}
\def \BPS {{\cal BP}_S}

\def \UM {{\cal UM}}
\def \UMX {{\cal UM}_X}

\def \UMS {{\cal UM}_S}

\def \nourm {the input does not yield a desired output}

\def \iff {if and only if }

\newcommand {\relabel}[1] {\label{#1}[\red{**its label}: \blue{\bf #1}]}\newcommand {\rebibitem}[1] {\bibitem{#1}[\red{**its label}: \blue{\bf #1}]} 

\def\relabel {\label} \def\rebibitem {\bibitem} 

\begin{document}
  
\newcommand{\resection}[1]
{\section{#1}\setcounter{equation}{0}}

\renewcommand{\theequation}{\thesection.\arabic{equation}}

\renewcommand{\labelenumi}{\rm(\roman{enumi})}

\baselineskip 0.6 cm

\title{From G-parking functions to B-parking functions\thanks{This paper was partially supported 
by NTU AcRF project (RP 3/16 DFM) of Singapore.}
}

\author
{ 
Fengming Dong\thanks{
Email: fengming.dong@nie.edu.sg.} \\
\small Mathematics and Mathematics Education\\
\small National Institute of Education\\
\small Nanyang Technological University, Singapore
\sgap \\
}

\date{}

\maketitle

\begin{abstract}
A matching $M$ in a multigraph $G=(V,E)$ is said to be  
uniquely restricted if $M$ is the only perfect matching 
in the subgraph of $G$ induced by $V(M)$ 
(i.e., the set of vertices saturated by $M$).
For any fixed vertex $x_0$ in $G$, 
there is a bijection from the set of
spanning trees of $G$ to the set of 
uniquely restricted matchings of size $|V|-1$ in $S(G)-x_0$,
where $S(G)$ is the bipartite graph obtained from $G$ 
by subdividing each edge in $G$.
Thus the notion ``uniquely restricted matchings
of a bipartite graph $H$ saturating all vertices in a 
partite set $X$" can be viewed as 
an extension of ``spanning trees in a connected graph".
Motivated by this observation, we extend the notion 
``G-parking functions" of a connected multigraph
to ``B-parking functions" 
$f:X\rightarrow \{-1,0,1,2,\cdots \}$
of a bipartite graph $H$ with a bipartition 
$(X,Y)$ 
and find a bijection $\psi$ from 
the set of uniquely restricted matchings 
of $H$ to the set of B-parking functions of $H$. 
We also show that for any 
uniquely restricted matching  $M$ in $H$ 
with $|M|=|X|$, if $f=\psi(M)$, 
then $\sum_{x\in X}f(x)$ is exactly the number of 
elements $y\in Y-V(M)$ which are not 
externally B-active with respect to $M$ in $H$,
where the new notion ``externally B-active members with 
respect to $M$ in $H$" is an extension of  
``externally active edges with respect to a spanning tree 
in a connected multigraph". 
\end{abstract}

\noindent {\bf MSC}: 05A19, 05B35 and 05C85

\noindent {\bf Keywords}:
graph, 
spanning tree,
parking function,
bijection

\resection{Introduction}

The notion of a parking function was introduced 
by Konheim and Weiss \cite{kon} in 1966. 
Suppose that there are $n$ drivers labeled $1,2,\cdots,n$
and $n$ parking spaces arranged in a line numbered $1,2,\cdots,n$.
Assume that these $n$ drivers enter the parking area in the order 
$1,2,\cdots,n$ and 
driver $i$ parks at space $j$, 
where $j$ is the minimum number 
with $f(i)\le j\le n$ such that space $j$ is unoccupied by the previous drivers
and $f(i)$ is the initial parking preference 
of driver $i$. 
If all drivers can park successfully by this rule, 
then $(f(1),f(2), \cdots,f(n))$ is called 
a {\it parking function} of length $n$.
Mathematically, a function $f:N_n\rightarrow N_n$,
where $N_n=\{1,2,\cdots,n\}$, 
is called {\it a parking function} 
if the inequality 
$|\{1\le i\le n: f(i)\le k\}|\ge k$ 
holds for each integer $k:1\le k\le n$.
For example, for $n=2$, 
$(f(1),f(2))=(1,1)$,  $(f(1),f(2))=(1,2)$
and $(f(1),f(2))=(2,1)$
are parking functions, but $(f(1),f(2))=(2,2)$ is not.
It can be shown easily that 
$f:N_n\rightarrow N_n$ is a parking function 
if and only if 
there is a permutation $\pi_1,\pi_2,\cdots,\pi_n$ 
of $N_n$ such that $f(\pi_j)\le j$ holds 
for all $j=1,2,\cdots,n$. 
Konheim and Weiss \cite{kon} proved that 
that the number of parking functions of length $n$ 
is equal to $(n+1)^{n-1}$, which is equal to the number 
of spanning trees of the complete graph $K_{n+1}$
(\cite{aig1998, cay1889}). 

The parking function and its various extensions 
have been studied by many researchers 
\cite{ma, che, ges1979, hai1994,  kos, kre1980, kru, lew1996, ma2016, pit2002, pos, prim, shi1986, shi1987, sim1994, sim1991, sta1996, sta1997, sta1998, ste, yan1997}.
One of the extensions,
due to  
Postnikov and Shapiro \cite{pos}, 
was from 
parking functions to G-parking functions 
for connected multigraphs without loops.

Before talking about  G-parking functions,
let's first introduce the notations of graphs used 
in this article. Unless stated otherwise, 
we always assume that 
\begin{enumerate}
\item $G=(V,E)$ is a connected multigraph
without loops, where $V=\{x_0,x_1,\cdots,x_n\}$
and $E=\{y_1,y_2,\cdots,y_m\}$.
For any non-empty subsets $V'$ of $V$ and 
$E'$ of $E$, let $G[V']$ and $G[E']$ be 
the subgraphs of $G$ induced by $V'$ and $E'$ respectively;

\item $H$ is a simple and bipartite graph with 
a bipartition $(X,Y)$, where $X=\{x_1,x_2,\cdots,x_n\}$
and $Y=\{y_1,y_2,\cdots,y_m\}$; and
\item $H_{G,x_0}$ is the special bipartite graph
$S(G)-x_0$ with a bipartition $(X,Y)$, 
where $x_0$ is a fixed vertex in $G$,
$S(G)$ is obtained from 
$G$ by subdividing each edge in $G$,
$X=V-\{x_0\}$ and $Y=E$.
An example of $H_{G,x_0}$ is shown in Figure~\ref{G-H}.
\end{enumerate}

Both graphs $G$ and $H$ have fixed weight functions
which are used for comparing edges in $G$ 
or elements of $Y$ in $H$.
The weight function for $G$ 
is an injective mapping $w:E\rightarrow \N_0$,
where $\N_0$ is the set of non-negative integers, 
while the weight function for $H$ is 
an injective mapping $w:Y\rightarrow \N_0$.
Thus the weight function $w$ of $G$ is also 
the weight function of $H_{G,x_0}$. 
The mapping $w$ is injective in order to 
distinguish $w(y_1)$ and $w(y_2)$
for any distinct elements $y_1$ and $y_2$.

For any  subsets $V_1$ and $V_2$ of $V$, 
let $E_G(V_1,V_2)$
denote the set of those edges in $G$ joining a vertex in $V_1$ and a vertex in $V_2$.
In particular, 
let $E_G(u,V_2)=E_G(\{u\},V_2)$ for any $u\in V$.
So $d_G(u)=|E_G(u,V)|$ is the degree of vertex $u$ in $G$.  
A function $f:V-\{x_0\}\rightarrow \N_0$ is called 
a {\it G-parking function} with respect to $x_0$ 
if for any non-empty subset $V'\subseteq V-\{x_0\}$,
there exists $u\in V'$ with $|E_G(u,V-V')|>f(u)$.
Let $\GP(G,x_0)$ denote the set of 
G-parking functions of $G$ 
with respect to $x_0$.

By Corollary~\ref{prop2-13}, 
which was due to Dhar~\cite{dhar},
a function $f:V-\{x_0\}\rightarrow \N_0$ 
belongs to $\GP(G,x_0)$ \iff 
there is an ordering $x_{\pi_1}, x_{\pi_2}, \cdots,
x_{\pi_n}$ 
of vertices in $V-\{x_0\}$  
such that $|E_G(x_{\pi_i},V-V_i)|>f(x_{\pi_i})$ 
holds for all $i=1,2,\cdots,n$,
where $V_i=\{x_{\pi_j}: i\le j\le n\}$.
Hence a function $f:N_n\rightarrow N_n$ is a 
parking function of length $n$ if and only if 
$f-1\in \GP(K_{n+1},0)$, 
where  $V(K_{n+1})=\{0,1,2,\cdots,n\}$.

The most interesting property on G-parking functions 
is the existence of bijections from 
the set of spanning trees of $G$, 
denoted by $\sett(G)$, to $\GP(G,x_0)$. 
Several such bijections have been obtained
(see~\cite{che} for example).

In this paper, we focus on presenting a new extension of 
G-parking functions.

A matching $M$ of a graph $G$ is said to be 
{\it uniquely restricted} (UR) if 
$M$ is the only perfect matching in $G[V(M)]$,
where $V(M)$ is the set of vertices saturated 
by edges in $M$.
Clearly, a matching $M$ of $G$ is a UR-matching
\iff $|E(C)|>2|E(C)\cap M|$ holds for every cycle 
$C$ in $G$,
where $E(C)$ is the set of edges on $C$.
The notion of UR-matchings 
was first introduced by 
Golumbic, Hirst, and Lewenstein \cite{gol},
originally motivated by the problem of determining a lower bound on the rank of a matrix having a specified zero/non-zero pattern.
They \cite{gol} showed that the problem of finding a UR-matching 
with the maximum cardinality in an
input graph is known to be NP-complete even for the special cases of split graphs and bipartite graphs.

For any $T\in \sett(G)$ with  
$E(T)=\{y_{\tau_i}:i=1,2,\cdots,n\}$, let 
$M_T$ denote the matching 
$\{x_{\pi_i}y_{\tau_i}: i=1,2,\cdots,n\}$ 
of $H_{G,x_0}$,
where $x_{\pi_i}$ is the end of edge $y_{\tau_i}$ in $G$ 
such that $y_{\tau_i}$ is contained in the unique path of $T$ 
connecting $x_0$ and $x_{\pi_i}$. 
An example of $T$ and $M_T$ is shown in Figure~\ref{G-H}.
Proposition~\ref{prop2-6} shows  
that the mapping $\lambda$ defined by 
$\lambda(T)=M_T$ is a bijection 
from  $\sett(G)$ 
to the set of UR-matchings of size $n$ ($=|V|-1$) in 
$H_{G,x_0}$.

\begin{figure}[htbp]
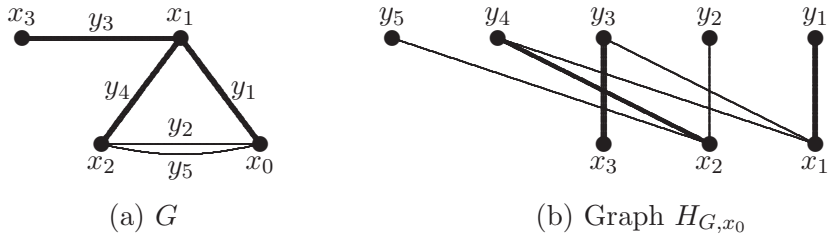

\centering

\begin{center}
\beginpicture 
\setcoordinatesystem units <2pt,2pt> 
\setplotarea x from -60 to 140, y from 0 to 20 
\plotsymbolspacing=.2pt 



\put {\cirs} at  -15 0  
\put {\cirs} at  15 0  
\put {\cirs} at  0 20   
\put {\cirs} at -30 20 

\plot 15.4 0 0.4 20 /  \plot 14.6 0 -0.4 20 /
\plot -14.6 0 0.4 20 /  \plot -15.4 0 -0.4 20 /
\plot -30 20.3 0.4 20.3 /  \plot -30 19.7 -0.4 19.7 /

\plot -15 0 15 0 0 20 /
\plot -0 20 -30 20  /
\plot -15 0 0 20 /

\put {$x_0$} at 15 -4
\put {$x_2$} at -15 -4
\put {$x_1$} at 0 24
\put {$x_3$} at -30 24
\put {$y_1$} at 12 10
\put {$y_4$} at -12 10
\put {$y_2$} at 0 3
\put {$y_5$} at 0 -5
\put {$y_3$} at -15 23



\put {\cirs} at  80 0
\put {\cirs} at  100 0
\put {\cirs} at  120 0

\put {$x_1$} at 120 -4
\put {$x_2$} at 100 -4
\put {$x_3$} at 80 -4

\put {\cirs} at  40 20
\put {\cirs} at  60 20
\put {\cirs} at  80 20
\put {\cirs} at  100 20
\put {\cirs} at  120 20

\linethickness= 8.500pt
\plot  120 0    80 20 /
\plot  120 0    60 20 /
\plot  100 0    60 20 /
\plot  100 0    40 20 /

\put {$y_1$} at 120 24
\put {$y_2$} at 100 24
\put {$y_3$} at 80 24
\put {$y_4$} at 60 24
\put {$y_5$} at 40 24

\linethickness= 4.500pt
\plot  100 0    100 20 /
\plot  120 0    120 20 /\plot 120.3 0 120.3 20 / \plot 119.7 0 119.7 20 /

\plot  80 0    80 20 /\plot 80.4 0 80.4 20 / \plot 79.65 0 79.65 20 /

\plot 100.5 0 60.5 20 / \plot 99.5 0 59.5 20 /

{\setquadratic 
\plot -15 0 0 -2 15 0 / 
}

\endpicture 
\end{center}

\vspace{-0.3 cm}

(a) $G$ \hspace{4.5 cm} (b) Graph $H_{G,x_0}$

\caption{$E(T)=\{y_1,y_3,y_4\}$ and 
$M_T=\{x_1y_1,x_3y_3,x_2y_4\}$}
\relabel{G-H}
\end{figure}

The above observation shows that 
the notion ``a spanning tree of a connected multigraph"
can be viewed as a special case of 
the notion ``a UR-matching of size $|X|$ in a 
bipartite graph $H$".
Motivated by this relation, we  extend 
the notion of G-parking functions of 
connected multigraphs to 
that of B-parking functions
of bipartite graphs. 

Let $\UM(H)$ be the set of UR-matchings of $H$.
For any $S\subseteq X$, 
let  $\UMS(H)$ be the set of those members $M$ 
of $\UM(H)$ with $V(M)\cap X=S$.
In particular, $\UMX(H)$ is the set of 
those members $M$ of $\UM(H)$
with $X\subseteq V(M)$. 
Thus $\UM(H)$ can be partitioned into subsets 
$\UMS(H)$ for all subsets $S$ of $X$.

A mapping 
$f:X\rightarrow \{-1\}\cup \N_0$ is called a 
{\it B-parking function} of $H$ at $X$
if for any non-empty subset $S$ of $X_{(f\ge 0)}$,
where $X_{(f\ge 0)}=\{x\in X: f(x)\ge 0\}$, 
there exists $x'\in S$ 
such that $x'$ has at least $f(x')+1$ neighbors of degree 1 
(i.e., leaves)
in the subgraph of $H$ induced by 
$\bigcup_{x\in S}N_H[x]$,
where $N_H(x)$ is the set of neighbors of $x$ in $H$ and 
$N_H[x]=\{x\}\cup N_H(x)$.
Let $\BP(H)$ be the family of B-parking functions of 
$H$ at $X$.
For any $S\subseteq X$, let $\BPS(H)$ 
be the set of those members $f\in \BP(H)$ 
with $X_{(f\ge 0)}=S$.
In particular, $\BPX(H)$ is the set of 
those members $f\in \BP(H)$  
with $f(x)\ge 0$ for all $x\in X$.
Thus $\BP(H)$  is also partitioned into
subsets $\BPS(H)$ for all subsets $S\subseteq X$.

In Section~\ref{sect2}, we give some basic properties 
on members in $\UMX(H)$ and members in $\BPX(H)$.
Proposition~\ref{prop2-12} shows that $\UMX(H)=\emptyset$ \iff $\BPX(H)=\emptyset$.
Proposition~\ref{prop2-6} shows that the members in $\sett(G)$ 
correspond to members in $\UMX(H_{G,x_0})$ 
and  Proposition~\ref{prop2-4} shows that 
$\GP(G,x_0)=\BPX(H_{G,x_0})$.

In Section~\ref{sect03},
we  design an algorithm, called Algorithm A, 
for any input $(H,Y')$, where $Y'\subseteq Y$. 
Whenever $\UMX(H[X\cup Y'])\ne \emptyset$, 
running this algorithm 
 outputs a permutation $\pi_1,\pi_2, \cdots, \pi_n$ 
of $1,2,\cdots,n$, an $n$-permutation $\tau_1,\cdots,\tau_n$ of 
$1,2,\cdots,m$ 
and subsets $D(x_{\pi_i})$ of $Y-Y'$ 
for $i=1,2,\cdots,n$.
In this case, 
the mapping $f: X\rightarrow \N_0$ defined by $f(x_{\pi_i})=|D(x_{\pi_i})|$ for $i=1,2,\cdots,n$ is a member in $\BPX(H)$.
This result yields a mapping 
$\psi_H$ from  $\UMX(H)$ to $\BPX(H)$.
The outputs $\pi_i, \tau_i$ and $D(x_{\pi_i})$ 
for $i=1,2,\cdots,n$ of running Algorithm A 
provide information for interpreting members in $\BPX(H)$.

In Section~\ref{sect04}, we show that 
the mapping  $\psi_H$ from  $\UMX(H)$ to $\BPX(H)$,
defined by $\psi_H(M)=f$,
is a bijection, 
where $f$ is the mapping from $X$ to $\N_0$ defined by  
$f(x_{\pi_i})=|D(x_{\pi_i})|$ for all $i=1,2,\cdots,n$,
and $\pi_i$ and $D(x_{\pi_i})$ are outputs of 
running Algorithm A with input $(H,V(M)\cap Y)$. 
Clearly, 
$\psi_{H[N[S]]}$ provides a bijection 
from $\UM_S(H)$ to $\BP_S(H)$ for every $S\subseteq X$,
where $N[S]=\bigcup_{x\in S}N_H[x]$.
Thus, 
there is a bijection from $\UM(H)$ to $\BP(H)$.
When $H$ is the graph $H_{G,x_0}$, 
$\psi_H$ is a bijection $\phi_G$ from $\sett(G)$ to 
$\GP(G,x_0)$ for any connected multigraph $G$, 
where $x_0\in V(G)$.

In Section~\ref{sect5}, we introduce 
the new notion 
``externally B-active members with respect to 
$M$ in $H$", where $M\in \UMX(H)$, 
defined in Page~\pageref{page-B-active},
which is an extension of  
``externally active edges with respect to a spanning tree 
$T$ in a connected multigraph" defined by Tutte~\cite{tut}.
For any $M\in \UMX(H)$, if $f=\psi_H(M)$, then 
$f(x_{\pi_i})$ is interpreted as the number of 
those $y\in N_H(x_{\pi_i})-\left (V(M)\cup  \bigcup_{s>i}N_H(x_{\pi_s})\right )$ 
which are not  
externally B-active with respect to $M$ in $H$,
implying that $\sum_{x_i\in X}f(x_i)$ is exactly the number of those vertices $y\in Y-V(M)$ which are not externally B-active with respect to $M$ in $H$.
This result implies that there exists a bijection 
$\phi_G$ from $\sett(G)$ to $\GP(G,x_0)$
such that for any $T\in \sett(G)$, 
if $f=\phi_G(T)$, then 
$\sum_{x\in V(G)-\{x_0\}}f(x)$ 
is exactly the number of those edges in $E(G)-E(T)$ 
which are not externally active with respect to $T$.

\resection{UR-matchings 
and B-parking functions\relabel{sect2}}

In this section, we  characterize UR-matchings 
and B-parking functions of a bipartite graph $H$. 
It is proved in Proposition~\ref{prop2-12} 
that  
$\UMX(H)=\emptyset$ \iff $\BPX(H)=\emptyset$. 
For the special bipartite graph $H_{G,x_0}$, 
Propositions~\ref{prop2-6} and~\ref{prop2-4} show 
that 
$\sett(G)$ and $\GP(G,x_0)$  
correspond to $\UMX(H_{G,x_0})$ and $\BPX(H_{G,x_0})$ respectively. 

\subsection{
UR-matchings in bipartite graphs
\relabel{sect2-1}}

By the definition of UR-matchings,  
the following statements are obviously equivalent 
for any matching $M$ in a multigraph $G$:
\begin{enumerate}
\item $M$ is a UR-matching of $G$;
\item 
$M$ is a UR-matching of the subgraph $G[V(M)]$;
\item $|E(C)|>2|M\cap E(C)|$ holds 
for any cycle $C$ in $G$.
\end{enumerate}

For UR-matchings in a bipartite graph,
another equivalent statement is given by 
Golumbic, Hirst and Hedetniemia~\cite{gol}.

\begin{theo}[\cite{gol}]\relabel{gol-theo}
$M\in \UMX(H)$
\iff
$M=\{x_{\pi_i}y_{\tau_i}: i=1,2,\cdots,n\}$
for a permutation $\pi_1,\pi_2,\cdots,\pi_n$ of $1,2,\cdots,n$ and 
an $n$-permutation 
$\tau_1,\tau_2,\cdots,\tau_n$ of $1,2,\cdots,m$ 
with $x_{\pi_i}y_{\tau_i}\in E(H)$ for all 
$i=1,2,\cdots,n$  
but $x_{\pi_j}y_{\tau_i}\notin E(H)$ for all 
$1\le i<j\le n$.
\end{theo}

Theorem~\ref{gol-theo}
can be stated equivalently as follows. 

\begin{cor}\relabel{cor2-11-0} 
For any $M\subseteq E(H)$ with $|M|=n$, 
$M\in \UMX(H)$
\iff 
$V(M)\cap Y=\{y_{\tau_i}: i=1,2,\cdots,n\}$ holds
for some $n$-permutation
$\tau_1,\tau_2,\cdots,\tau_n$ of $1,2,\cdots,m$ 
such that  
$y_{\tau_i}$ is a leaf in the subgraph 
$H-\bigcup_{1\le s<i}N_H[y_{\tau_s}]$
for all $i=1,2,\cdots,n$.
\end{cor}

Corollary~\ref{cor2-11-0} implies 
a necessary condition for 
$\UMX(H)$ to be non-empty.
Let $L(H)$ denote the set of leaves in $H$.

\begin{cor}\relabel{cor2-11-1} 
If $\UMX(H)\ne \emptyset$,
then $L(H_i)\cap Y\ne \emptyset $ for 
each component $H_i$ of $H$.
\end{cor}

But Corollary~\ref{cor2-11-1} 
is not true  for a non-bipartite graph
which contains  perfect UR-matchings. 
An example from \cite{gol} 
is shown in Figure~\ref{f10},
where the graph is non-bipartite and 
has a perfect UR-matching $\{e_1,e_2,e_3\}$. 
But it does not have any leaf. 

\begin{figure}[htbp]
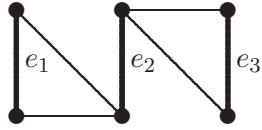

\centering

\beginpicture 
\setcoordinatesystem units <2pt,2pt> 
\setplotarea x from -110 to 100, y from 0 to 20 
\plotsymbolspacing=.2pt

\put {\cirs} at  -20 0
\put {\cirs} at  0 0
\put {\cirs} at  20 0

\put {\cirs} at  -20 20
\put {\cirs} at  0 20
\put {\cirs} at  20 20

\plot -20.3 0 -20.3 20 /
\plot -19.7 0 -19.7 20 /

\plot 0.3 0 0.3 20 /
\plot -0.3 0 -0.3 20 /

\plot 20.3 0 20.3 20 /
\plot 19.7 0 19.7 20 /

\plot 0 0 -20 0 -20 20 /
\plot -20 20 0 0 0 20 /
\plot 20 0 20 20 0 20 20 0 /

\put {$e_1$} at  -16 10
\put {$e_2$} at  4 10
\put {$e_3$} at  24 10

\endpicture 


\caption{A non-bipartite graph with a 
perfect UR-matching $\{e_1,e_2,e_3\}$}
\relabel{f10}
\end{figure}

Hall's Theorem~\cite{hall1935} on bipartite graphs
is an important result of characterizing bipartite
graphs with matchings saturating all vertices in one 
partite set. 
By Theorem~\ref{gol-theo}, 
we can get a characterization 
for $\UMX(H)$ to be non-empty
in terms of the sizes of sets $N_H(S)$, where $S\subseteq X$.

\begin{cor}\relabel{cor2-12}
$\UMX(H)\ne \emptyset$
\iff
there exists a permutation $\pi_1,\pi_2,\cdots,\pi_n$ 
of $1,2,\cdots,n$ such that 
$$
|N_H(X_1)|>|N_H(X_2)|>\cdots>|N_H(X_n)|>0,
$$
where 
$X_i=\{x_{\pi_j}:i\le j\le n\}$
and $N_H(X_i)=\bigcup_{x\in X_i}N_H(x)$.
\end{cor}

By Corollary~\ref{cor2-12} or Theorem~\ref{gol-theo}, 
if $\UMX(H)\ne \emptyset$, then 
$H$ contains at least one leaf $y'\in Y$ in $H$.
We are now going to show that 
if $\UMX(H)\ne \emptyset$,
then each leaf $y'\in Y$ of $H$
is contained in $V(M)$ for some 
$M\in \UMX(H)$.

\begin{pro}\relabel{prop2-01}
Assume that $y'\in Y$ is a leaf of $H$ 
with $N_H(y')=\{x'\}$.
Let $X'=X-\{x'\}$, $H'=H-y'$ and $H''=H-\{x',y'\}$. 
The following statement hold:
\begin{enumerate}
\item  for any $M\in \UMX(H)$,
if $y'\notin V(M)$, then $M\in \UMX(H')$;
otherwise,  
$M-\{x'y'\}\in \UM_{X'}(H'')$;

\item if $\UMX(H)\ne \emptyset$, then 
$y'\in V(M)$ for some $M\in \UMX(H)$.
\end{enumerate}  
\end{pro}

\proof 
(i) follows from Theorem~\ref{gol-theo} directly.

(ii) 
Assume that $M\in  \UMX(H)$ with $y'\notin V(M)$.
By Theorem~\ref{gol-theo},
there exist a permutation $\pi_1,\pi_2,\cdots,\pi_n$ 
of $1,2,\cdots,n$ and an $n$-permutation 
$\tau_1,\tau_2,\cdots,\tau_n$ of 
$1,2,\cdots,m$  
such that $M=\{x_{\pi_i}y_{\tau_i}: i=1,2,\cdots,n\}$  
and $x_{\pi_i}y_{\tau_j}\notin E(H)$ for all 
$1\le j<i\le n$.

Assume that $y'=y_q$ and $x'=x_{\pi_k}$.
Then $\tau_i\ne q$ for all $i=1,2,\cdots,n$.
Let $\gamma_k=q$ and 
$\gamma_i=\tau_i$ for all $i$ with $1\le i\le n$ 
and $i\ne k$.
Then $\pi_1,\pi_2,\cdots,\pi_n$ is a permutation 
of $1,2,\cdots,n$ and 
$\gamma_1,\gamma_2,\cdots,\gamma_n$ is 
an $n$-permutation of $1,2,\cdots,m$  
such that 
$x_{\pi_i}y_{\gamma_i}\in E(H)$ for all $i=1,2,\cdots,n$
but 
$x_{\pi_i}y_{\gamma_j}\notin E(H)$ for all 
$1\le j<i\le n$.
By Theorem~\ref{gol-theo}, 
$M'=\{x_{\pi_i}y_{\gamma_i}, i=1,2,\cdots,n\}$
is a member in $\UMX(H)$ 
with $y'=y_{q}=y_{\gamma_k}\in V(M')$.

Hence (ii) holds.
\proofend

\subsection{B-parking functions}

A characterization of B-parking functions is given below.

\begin{pro}\relabel{prop2-11}
For any mapping $f:X\rightarrow \N_0$, 
$f\in \BPX(H)$ \iff 
there is a permutation $\pi_1,\pi_2,\cdots,\pi_n$ 
of $1,2,\cdots,n$ such that for each $i=1,2,\cdots,n$,
$x_{\pi_i}$ has at least $f(x_{\pi_i})+1$ neighbors 
which are leaves in the subgraph of $H$ 
induced by $\bigcup_{i\le j\le n}N[x_{\pi_j}]$.
\end{pro}

\proof 
($\Rightarrow $) 
Assume that $f\in \BPX(H)$.
By the definition of B-parking functions, 
there exists a vertex $x_{\pi_1}\in X$
such that $|N_H(x_{\pi_1})\cap L(H)|\ge f(x_{\pi_1})+1$.
 
Assume that $\pi_1,\pi_2,\cdots,\pi_s$ 
is a $s$-permutation of $1,2,\cdots,n$, where $1\le s<n$,
such that for all $i=1,2,\cdots,s$,
$|N_H(x_{\pi_i})\cap L(H[N[X_i]])|\ge f(x_{\pi_i})+1$,
where $X_i=X-\{x_{\pi_r}: 1\le r<i\}$.
By the definition of B-parking functions again, 
there exists a vertex, denoted by $x_{\pi_{s+1}}$,
in $X_{s+1}$ such that 
$|N_H(x_{\pi_{s+1}})\cap L(H[N[X_{s+1}]])|
\ge f(x_{\pi_{s+1}})+1$.
Repeating this process, 
a permutation $\pi_1,\pi_2,\cdots,\pi_n$ of $N_n$ 
can be obtained 
such that  
$|N_H(x_{\pi_i})\cap L(H[N[X_i]])|\ge f(x_{\pi_i})+1$
for all $i=1,2,\cdots,n$.
Observe that $X_i$ is the set 
$\{x_{\pi_{r}}:i\le r\le n\}$. 
Thus the necessity holds.

($\Leftarrow$) 
Now assume that $\pi_1,\pi_2,\cdots,\pi_n$ 
is a permutation 
of $1,2,\cdots,n$ such that for $i=1,2,\cdots,n$,
$|N_H(x_{\pi_i})\cap L(H[N[X_i]])|\ge f(x_{\pi_i})+1$
holds, 
where  $X_i=\{x_{\pi_{r}}:i\le r\le n\}$. 
Let $X'$ be an arbitrary non-empty subset of $X$
and $s$ be the minimum integer in $N_n$ such that 
$x_{\pi_s}\in X'$. 
By assumption, 
$x_{\pi_s}$ has at least $f(x_{\pi_s})+1$ neighbors 
which are leaves in $H[N[X_s]]$.  
Observe that $X'\subseteq X_s=\{x_{\pi_r}: s\le r\le n\}$,
implying that for any $y\in N_H(x_{\pi_s})$, 
if $y\in L(H[N[X_s]])$,
then $y\in L(H[N[X']])$.
Thus $|N_H(x_{\pi_s})\cap L(H[N[X']])\ge f(x_{\pi_s})+1$.
Hence $f\in \BPX(H)$.
\proofend

By Proposition~\ref{prop2-11}, one can prove the following characterization for members in $\BPX(H)$ 
by acyclic orientations of $H$.

\begin{cor}\relabel{cor2-02-0}
For any $f:X\rightarrow \N_0$, 
$f\in \BPX(H)$ if and only if 
there exists an acyclic orientation $D$ of $H$ 
such that $od_D(y_j) =1$ 
holds for all $j=1, 2, \cdots, m$ and 
$f(x_i) <id_D(x_i)$ holds for all $i =1, 2, \cdots, n$, 
where $od_D(y_j)$ and $id_D(x_i)$ 
are respectively the out-degree of $y_j$ 
and the in-degree of $x_i$ in $D$.
\end{cor}

Let $f$ be a mapping from $X$ to $\N_0$.
For any $X'\subseteq X$ and $x'\in X$, 
let $f|_{X'}$ be 
the restriction of $f$ to the set $X'$
and let $f_{(x'\downarrow 1)}$ be the mapping defined by  
$f_{(x'\downarrow 1)}(x')=f(x')-1$ 
and $f_{(x'\downarrow 1)}(x)=f(x)$ for all $x\in X-\{x'\}$.
By Proposition~\ref{prop2-11}, 
we have the following result.

\begin{cor}\relabel{cor2-02}
Assume that $y'\in Y\cap L(H)$ and 
$N_H(y')=\{x'\}$. 
For any mapping $f$ from $X$ to $\N_0$,
the following statements hold:
\begin{enumerate}
\item if $f(x)=0$ for all $x\in X$,
then $\BPX(H)\ne \emptyset$ \iff 
$f\in \BPX(H)$;

\item  $f_{(x'\downarrow 1)}\in \BPX(H-y')$ 
\iff $f\in \BPX(H)$ and $f(x')\ge 1$;  

\item if $f(x')=0$,
then $f|_{X-\{x'\}}\in \BP_{X-\{x'\}}(H-x')$ \iff $f\in \BPX(H)$.
\end{enumerate}  
\end{cor}

By applying Propositions~\ref{prop2-01} and~\ref{prop2-11},
Theorem~\ref{gol-theo} and Corollary~\ref{cor2-02},
it can be shown that 
$\UMX(H)\ne \emptyset$ 
\iff $\BPX(H)\ne \emptyset$.

\begin{pro}\relabel{prop2-12}
The following statements are equivalent:
\begin{enumerate}
\item $L(H)\cap Y\ne \emptyset$ and 
for each $y\in L(H)\cap Y$, 
$y\in V(M)$ for some $M\in \UMX(H)$;
\item $\UMX(H)\ne \emptyset$;
\item
there exist a permutation $\pi_1, \pi_2, \cdots, \pi_n$ 
of $1,2,\cdots,n$ and 
an $n$-permutation $\tau_1,\tau_2,\cdots,\tau_n$ 
of $1,2,\cdots,m$ such that 
$M=\{x_{\pi_i}y_{\tau_i}:i=1,2,\cdots,n\}$ 
and $x_{\pi_i}y_{\pi_j}\notin E(H)$ for all 
$1\le j<i\le n$;
\item
$f\in \BPX(H)$, where $f$ is the mapping 
defined by $f(x)=0$ for all $x\in X$;
\item $\BPX(H)\ne \emptyset$.
\end{enumerate}
\end{pro}

\proof
Observe that (i) $\Leftrightarrow $ (ii), 
(ii) $\Leftrightarrow $ (iii),
(iii) $\Leftrightarrow $ (iv) and 
(iv) $\Leftrightarrow $ (v) follow from 
Proposition~\ref{prop2-01} (ii), 
Theorem~\ref{gol-theo}, 
Proposition~\ref{prop2-11}
and Corollary~\ref{cor2-02} (i) respectively.
\proofend

\subsection{$\UMX(H_{G,x_0})$ and $\BPX(H_{G,x_0})$
\relabel{sect2-3} }

We focus on the special bipartite graph $H_{G,x_0}$
in this subsection.
Note that $H_{G,x_0}$ has a bipartition $(X,Y)$,
where $X=V-\{x_0\}$ and $Y=E$.
Each vertex of $Y$ is of degree $1$ or $2$ in $H_{G,x_0}$.
As $G$ is connected, 
$L(H_t)\cap Y\ne \emptyset$ for 
each component $H_t$ of $H_{G,x_0}$.
Also note that $y_i$ and $y_j$ are parallel edges in 
$G$ \iff $y_i$ and $y_j$ have the same set 
of neighbors in $H_{G,x_0}$.
An example of $H_{G,x_0}$ is shown in Figure~\ref{f5-3}.

In this subsection, 
we will show that 
there is a bijection from $\sett(G)$ to $\UMX(H_{G,x_0})$
and $\GP(G,x_0)=\BPX(H_{G,x_0})$ holds.

\begin{lem}\relabel{le2-3-1}
If $G_0$ is a disconnected multigraph, then $\UMX(H_{G_0,x_0})=\emptyset$.
\end{lem}

\proof Assume that $G_0$ is disconnected. 
Then some component  of $H_{G_0,x_0}$ 
does not have leaves.
By Corollary~\ref{cor2-11-1}, 
$\UMX(H_{G_0,x_0})=\emptyset$.
\proofend

By Lemma~\ref{le2-3-1}, we need only to 
consider connected multigraphs. 
Let
$T\in \sett(G)$. 
Without loss of generality,
 assume that $E(T)=\{y_i:1\le i\le n\}$.
Recall that 
$M_T$ denotes the matching 
$\{x_{\epsilon_i}y_{i}: i=1,2,\cdots,n\}$
of $H_{G,x_0}$,
where $x_{\epsilon_i}$ is the end of edge $y_{i}$ in $G$ 
such that $y_{i}$ is contained in the unique path 
in $T$ from $x_0$ to $x_{\epsilon_i}$.
By the definition of $M_T$, 
$M_T$ is characterized by the following lemma.

\begin{lem}\relabel{le2-3}
$M_T=\{x_{\pi_i}y_{\tau_i}: i=1,2,\cdots,n\}$
\iff   
$\pi_1,\pi_2,\cdots,\pi_n$  is a permutation 
of $1,2,\cdots,n$ such that 
each $y_{\tau_i}$ is an edge in $E_T(V_i,V-V_i)$
incident with $x_{\pi_i}$,
where 
$V_i=\{x_0\}\cup \{x_{\pi _j}: 1\le j<i\}$
for $i=1,2,\cdots,n$.
\end{lem}

\begin{figure}[htbp]
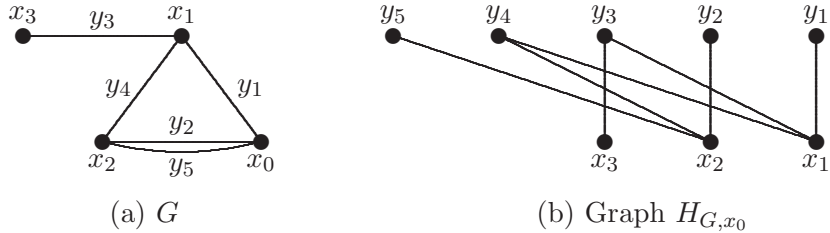

\centering

\begin{center}
\beginpicture 
\setcoordinatesystem units <2pt,2pt> 
\setplotarea x from -60 to 140, y from 0 to 20 
\plotsymbolspacing=.2pt 



\put {\cirs} at  -15 0
\put {\cirs} at  15 0
\put {\cirs} at  0 20 
\put {\cirs} at -30 20 

\plot -15 0 15 0 0 20 /
\plot -0 20 -30 20  /
\plot -15 0 0 20 /

\put {$x_0$} at 15 -4
\put {$x_2$} at -15 -4
\put {$x_1$} at 0 24
\put {$x_3$} at -30 24
\put {$y_1$} at 13 10
\put {$y_4$} at -12 10
\put {$y_2$} at 0 3
\put {$y_5$} at 0 -5
\put {$y_3$} at -15 23



\put {\cirs} at  80 0
\put {\cirs} at  100 0
\put {\cirs} at  120 0

\put {$x_1$} at 120 -4
\put {$x_2$} at 100 -4
\put {$x_3$} at 80 -4

\put {\cirs} at  40 20
\put {\cirs} at  60 20
\put {\cirs} at  80 20
\put {\cirs} at  100 20
\put {\cirs} at  120 20

\linethickness= 8.500pt

\plot  120 0    80 20 /
\plot  120 0    60 20 /
\plot  100 0    60 20 /
\plot  100 0    40 20 /

\put {$y_1$} at 120 24
\put {$y_2$} at 100 24
\put {$y_3$} at 80 24
\put {$y_4$} at 60 24
\put {$y_5$} at 40 24

\linethickness= 4.500pt
\plot  80 0    80 20 /
\plot  100 0    100 20 /
\plot  120 0    120 20 /

{\setquadratic 
\plot -15 0 0 -2 15 0 / 
}

\endpicture 
\end{center}

\vspace{-0.3 cm}

(a) $G$ \hspace{4.5 cm} (b) Graph $H_{G,x_0}$

\caption{Graphs $G$ and $H_{G,x_0}$}
\relabel{f5-3}
\end{figure}

\begin{pro}\relabel{prop2-6}
The mapping $\lambda: \sett(G)\rightarrow \UMX(H_{G,x_0})$
defined by $\lambda(T)=M_T$ is a bijection. 
\end{pro}

\proof Clearly, if $T_1$ and $T_2$ are 
distinct members in $\sett(G)$, then $M_{T_1}\ne M_{T_2}$.
Thus, it suffices to prove the following statements:
\begin{enumerate}
\item For any $T\in \sett(G)$,  $M_T\in \UMX(H_{G,x_0})$;
\item For any $T\in \sett(G)$,  
$M_T$ is the only member in $\UMX(H_{G,x_0})$
with $V(M_T)\cap Y=E(T)$;
\item For any $M\in \UMX(H_{G,x_0})$,
$V(M)\cap Y=E(T)$ holds for some $T\in \sett(G)$. 
\end{enumerate}

(i) 
Let $T\in \sett(G)$. 
By Lemma~\ref{le2-3},
$M_T=\{x_{\pi_i}y_{\tau_i}: i=1,2,\cdots,n\}$, 
where $\pi_1,\pi_2,\cdots,\pi_n$ is some permutation 
of $1,2,\cdots,n$ such that 
$y_{\tau_i}$ is an edge in $E_T(V_i,V-V_i)$
with $x_{\pi_i}$ as one end 
and $V_i=\{x_0\}\cup \{x_{\pi _j}: 1\le j<i\}$
for all for $i=1,2,\cdots,n$.
Thus $y_{\tau_i}$ is a leaf in 
$H_{G,x_0}-\bigcup_{1\le j<i}N_{H_{G,x_0}}(y_{\tau_i})$
for all $i=1,2,\cdots,n$.
As $M_T\subseteq E(H_{G,x_0})$, 
by Corollary~\ref{cor2-11-0},
$M_T\in \UMX(H_{G,x_0})$.
Thus (i) holds.

(ii) For any $T\in \sett(G)$, 
by definition, $H_{T,x_0}$ is exactly the subgraph 
$H_{G,x_0}[X\cup E(T)]$.
Since $H_{T,x_0}=S(T)-x_0$ has no cycles, 
$H_{G,x_0}[X\cup E(T)]$ has no 
cycles, implying that $H_{G,x_0}[X\cup E(T)]$
cannot have two distinct perfect matchings. 
As each member $M\in \UMX(H)$ with $V(M)\cap X=E(T)$
is a perfect matching of $H_{G,x_0}[X\cup E(T)]$,
(ii) holds. 

(iii) Let $M\in \UMX(H_{G,x_0})$ and $Y'=V(M)\cap Y$.
By  Corollary~\ref{cor2-11-0},
there is an $n$-permutation $\tau_1,\tau_2,\cdots, \tau_n$ 
of $1,2,\cdots,m$ such that for all $i=1,2,\cdots,n$,
$y_{\tau_i}\in Y'$ and $y_{\tau_i}$ is 
incident with a unique vertex $x_{\pi _i}$
in the subgraph 
$H_{G,x_0}-\bigcup_{1\le j<i}N_{H_{G,x_0}}(y_{\tau_i})$,
implying that $y_{\tau_i}\in E_G(V_i,V-V_i)$
with $x_{\pi_i}$ as one end,
where $V_i=\{x_0\}\cup \{x_{\pi_j}: 1\le j<i\}$.
Thus,  $G[Y']$ is a tree.
Hence (iii) holds. 
\proofend

Proposition~\ref{prop2-6} shows that 
the notion of 
UR-matchings in bipartite graphs is an 
extension of that of spanning trees 
in connected multigraphs.

\begin{pro}\relabel{prop2-4}
For any mapping $f: X\rightarrow \N_0$, 
$f\in \GP(G,x_0)$ \iff $f\in \BPX(H_{G,x_0})$. 
\end{pro}

\proof 
Consider the following statements:
\begin{enumerate}
\item $f\in \GP(G,x_0)$;
\item for any non-empty subset $V'$ of $X$, 
there exists $x_j\in V'$ with
$|E_G(x_j,V-V')|>f(x_j)$;
\item 
for any non-empty subset $V'$ of $X$, 
there exists $x_j\in V'$ 
such that $x_j$ has at least $f(x_j)$ 
neighbors which are leaves
in the subgraph of $H_{G,x_0}$ induced by 
$\bigcup_{x_i\in V'} N_{H_{G,x_0}}(x_i)$;
\item 
$f\in \BPX(H_{G,x_0})$. 
\end{enumerate}
(i) $\Leftrightarrow $ (ii) 
and (iii) $\Leftrightarrow $ (iv)  
follow from the definitions
of $\GP(G,x_0)$ and $\BPX(H_{G,x_0})$ respectively.
(ii) $\Leftrightarrow $ (iii)
follows from the fact that $y\in E_G(x_j,V-V')$
\iff $y$ is a vertex in $H_{G,x_0}$ adjacent to $x_j$
and is also a leaf in the subgraph of $H_{G,x_0}$ 
induced by $\bigcup_{x_i\in V'}N[x_i]$.
Hence the result holds.
\proofend

A characterization on G-parking functions
follows directly from Proposition~\ref{prop2-11} and 
Proposition~\ref{prop2-4}.
It was first obtained by Dhar~\cite{dhar}.

\begin{cor}[Dhar~\cite{dhar}]\relabel{prop2-13}
For any $f: V-\{x_0\}\rightarrow \N_0$, 
$f\in \GP(G,x_0)$ \iff 
 there is a permutation $\pi_1, \pi_2, \cdots,\pi_n$ 
of $1,2,\cdots,n$
such that  
$|E_G(x_{\pi_i},V-V_i)|>f(x_{\pi_i})$ holds
for each $i=1,2,\cdots,n$,
where $V_i=\{x_{\pi_j}: i\le j\le n\}$.
\end{cor}

Applying the notion of acyclic orientations of $G$,
Corollary~\ref{prop2-13} can be equivalently stated 
as follows.

\begin{cor}\relabel{G-ori}
For any $f: V-\{x_0\}\rightarrow \N_0$, 
$f\in \GP(G,x_0)$ \iff 
there exists an acyclic orientation $D$ of $G$ 
with $x_0$ as its unique source 
such that $f(x_i)<id_D(x_i)$ holds 
for all $i=1,2,\cdots,n$,
where $id_D(x_i)$ is the in-degree of $x_i$ in $D$. 
\end{cor}

\resection{An algorithm
\relabel{sect03}}

In this section, 
we design an algorithm, called  {\it Algorithm A},
mainly for the purpose of 
producing a member $f$ in $\BPX(H)$  
for any $Y'\subseteq Y$ 
with $\UMX(H[X\cup Y'])\ne \emptyset$, 
as stated in Proposition~\ref{prop2-44}.
By this result, we are able to define a mapping $\psi_H$ 
from $\UMX(H)$ to $\BPX(H)$ which is shown to be a bijection 
in Theorem~\ref{B-bi}.
The outputs of this algorithm are also applied in Section~\ref{sect5} to 
interpret the member $f\in \BPX(H)$ which corresponds to 
any given $M\in \UMX(H)$ under the mapping $\psi_H$. 

\subsection{Algorithm A\relabel{sect3-1}}

The weight function $w:Y\rightarrow \N_0$ of $H$ 
is needed for running Algorithm A.
In order to distinguish members in $Y$, 
we assume that $w$ is injective and so 
$w(y_1)\ne w(y_2)$ holds for 
any two different members $y_1,y_2\in Y$.
The input for 
Algorithm A below is an order pair $(H,Y')$, 
where $Y'\subseteq Y$.

\noindent {\bf Algorithm A} $(H, Y')$:
\begin{enumerate}
\item[A1:] Input $H$ with 
a bipartition $(X,Y)$ and a subset $Y'$ of $Y$;
\item[A2:] Set $i := 1$, 
$I:=X$, 
$D(x):= \emptyset$ and 
$F(x):= N_H(x)$  for all $x\in X$;
\item[A3:] Set 
$$
L_I:=\{y\in \bigcup_{x\in I}F(x): y \mbox{ is a leaf in } H_I\},
$$
where $H_I$ is the subgraph of $H$ induced by 
$I\cup \left (\bigcup_{x\in I}F(x)\right )$.
If $L_I=\emptyset$,
then output the message ``\nourm" and stop;
\item[A4:] If $L_I\ne \emptyset$,
determine the member $y'$ in $L_I$ 
with $w(y')<w(y)$ for all $y\in L_I-\{y'\}$
and the unique member $x'\in N_H(y')$;

\item[A5:] If $y'\notin Y'$, then 
set $F(x'):=F(x')-\{y'\}$, 
$D(x'):=D(x')\cup\{y'\}$ and go back to Step A3;

\item[A6:] If $y'\in Y'$,
determine the unique number $\pi_i\in \{1,2,\cdots,n\}$ 
and the unique number  $\tau_i\in \{1,2,\cdots,m\}$
such that $x_{\pi_i}=x'$ and $y_{\tau_i}=y'$; 

\item[A7:]
Set  $I:=I-\{x'\}$.
If $|I|>0$, set $i:=i+1$ and go back to Step A3; 

\item[A8:] Output $\pi_i,\tau_i$ and $D(x_{\pi_i})$
for all $i=1,2,\cdots,n$ and stop.
\end{enumerate}

Running Algorithm A has two possible outcomes.
Let $\sigma(H,Y')=0$ if running Algorithm A with inputs 
$(H,Y')$
stops with the message ``\nourm",
and let $\sigma(H,Y')=1$ otherwise.
In the case $\sigma(H,Y')=1$,
running Algorithm A outputs numbers $\pi_i$,
$\tau_i$ and a subset $D(x_{\pi_i})$ of $Y-Y'$ 
for $i=1,2,\cdots,n$, 
where $\pi_1,\pi_2,\cdots,\pi_n$ 
is a permutation of $1,2,\cdots,n$ 
and $\tau_1,\tau_2,\cdots,\tau_n$ 
is an $n$-permutation of $1,2,\cdots,m$. 
In this case,  $\pi_i$,  
$\tau_i$ and  $D(x_{\pi_i})$  
are rigorously written as 
$\pi_i(H,Y')$, $\tau_i(H,Y')$
and $D(H,Y',x_{{\pi_i}})$. 

\begin{table}[htbp]
\begin{center}
\begin{tabular}{c|c|c|c|c}
& $i=1$ &$i=2$ & $i=3$& $i=4$ \\ \hline
$\pi_i(H_1,Y_1)$ & 4 & 3& 2& 1\\ \hline
$\tau_i(H_1,Y_1)$ & 5 & 6&1 & 2\\ \hline
$D(H_1,Y_1,x_{\pi_i})$ & $\emptyset$ & $\emptyset$ & $\emptyset$ & $\emptyset$ \\ 
\end{tabular}
\end{center}
\caption{
$\pi_i(H_1,Y_1)$, $\tau_i(H_1,Y_1)$
and $D(H_1,Y_1,x_{\pi_i})$,
where $Y_1=\{y_1,y_2,y_5,y_6\}$}

\relabel{table1}
\end{table}

Let's consider some examples. 
Let $H_1$ and $H_2$ be bipartite graphs given 
in Figure~\ref{f12} with $w(y_i)=i$.
It is not difficult to verify that 
$\sigma(H_2,Y')=0$ for all subsets $Y'$
of $\{y_1,y_2,y_3,y_4,y_5\}$.
For graph $H_1$, we have 
$\sigma(H_1,Y')=0$ if $Y'=\{y_1,y_2,y_3,y_4\}$.
But 
$\sigma(H_1,Y_i)=1$ for $i=1,2$,
where $Y_1=\{y_1,y_2,y_5,y_6\}$
and $Y_2=\{y_3,y_4,y_5,y_6\}$,
and the outputs are shown in Tables~\ref{table1}
and~\ref{table2} respectively. 

\begin{table}[htbp]
\begin{center}
\begin{tabular}{c|c|c|c|c}
& $i=1$ &$i=2$ & $i=3$& $i=4$ \\ \hline
$\pi_i(H_1,Y_2)$ & 4 & 3& 1& 2\\ \hline
$\tau_i(H_1,Y_2)$ & 5 & 6&4 & 3\\ \hline
$D(H_1,Y_2,x_{\pi_i})$ & $\emptyset$ & $\emptyset$ & 
 $\{y_2\}$ & $\{y_1\}$ \\ 
\end{tabular}
\end{center}

\caption{
$\pi_i(H_1,Y_2)$, $\tau_i(H_1,Y_2)$
and $D(H_1,Y_2,x_{\pi_i})$,
where $Y_2=\{y_3,y_4,y_5,y_6\}$}
\relabel{table2}
\end{table}

\begin{figure}[htbp]
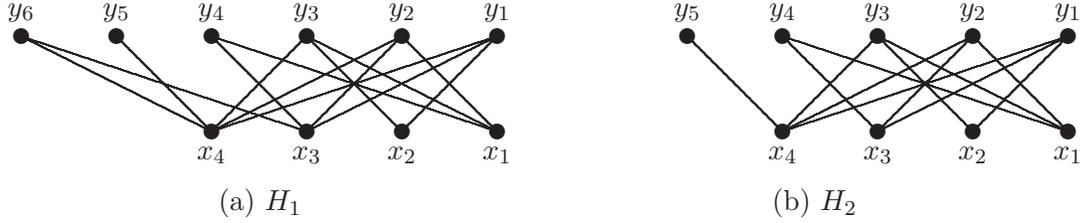

\centering

\beginpicture 
\setcoordinatesystem units <1.8pt,1.8pt> 
\setplotarea x from -130 to 100, y from -2 to 20 
\plotsymbolspacing=.18pt

\put {\cirs} at  -120 20  \put {$y_6$} at  -120 25
\put {\cirs} at  -100 20 \put {$y_5$} at  -100 25
\put {\cirs} at  -80 20 \put {$y_4$} at  -80 25
\put {\cirs} at  -60 20 \put {$y_3$} at  -60 25
\put {\cirs} at  -40 20 \put {$y_2$} at  -40 25
\put {\cirs} at  -20 20 \put {$y_1$} at  -20 25

\put {\cirs} at  -80 0  \put {$x_4$} at  -80 -5
\put {\cirs} at  -60 0 \put {$x_3$} at  -60 -5
\put {\cirs} at  -40 0 \put {$x_2$} at  -40 -5
\put {\cirs} at  -20 0 \put {$x_1$} at  -20 -5

\plot -100 20 -80 0 -120 20  -60 0 /  
\plot -20 20 -80 0 / 
\plot -40 20 -80 0 / 
\plot -60 20 -80 0 / 

\plot -20 20 -60 0 / 
\plot -40 20 -60 0 / 
\plot -80 20 -60 0 / 

\plot -20 20 -40 0 / 
\plot -60 20 -40 0 / 

\plot -40 20 -20 0 / 
\plot -60 20 -20 0 / 
\plot -80 20 -20 0 / 



\put {\cirs} at  20 20 \put {$y_5$} at  20 25
\put {\cirs} at  40 20 \put {$y_4$} at  40 25
\put {\cirs} at  60 20 \put {$y_3$} at  60 25
\put {\cirs} at  80 20 \put {$y_2$} at  80 25
\put {\cirs} at  100 20 \put {$y_1$} at  100 25

\put {\cirs} at  40 0  \put {$x_4$} at  40 -5
\put {\cirs} at  60 0 \put {$x_3$} at  60 -5
\put {\cirs} at  80 0 \put {$x_2$} at  80 -5
\put {\cirs} at  100 0 \put {$x_1$} at  100 -5

\plot 20 20 40 0 /  
\plot 100 20 40 0 / 
\plot 80 20 40 0 / 
\plot 60 20 40 0 / 

\plot 100 20 60 0 / 
\plot 80 20 60 0 / 
\plot 40 20 60 0 / 

\plot 100 20 80 0 / 
\plot 60 20 80 0 / 

\plot 80 20 100 0 / 
\plot 60 20 100 0 / 
\plot 40 20 100 0 / 

\endpicture 

\vspace{0.3 cm}

(a) $H_1$ \hspace{6 cm} (b) $H_2$

\caption{Bipartite graphs $H_1$ and $H_2$}
\relabel{f12}
\end{figure}

\subsection{When does the case ``$\sigma(H,Y')=1$" happen
\relabel{sect3-2}}

In this subsection, we shall know when the case ``$\sigma(H,Y')=1$" happens,
and how the outputs $\pi_i,\tau_i$ and $D(x_{\pi_i})$ are determined when it happens.

If $L(H)\cap Y=\emptyset$, 
then $\sigma(H,Y')=0$ clearly. 
If $L(H)\cap Y\ne \emptyset$, 
we have the following observations
from Algorithm A.

\begin{lem}\relabel{le2-2}
Assume that $L(H)\cap Y\ne \emptyset$ and 
$y'$ is the member in $L(H)\cap Y$ such that 
$w(y')$ is the minimum. 
Let $Y'\subseteq Y$, $Y''=Y'-\{y'\}$,
$H'=H-y'$ and $H''=H-\{x',y'\}$,
where $x'$ is the only member in $N_H(y')$.
The following observations follow from Algorithm A:
\begin{enumerate}
\item  if $y'\notin Y'$, then $\sigma(H,Y')=\sigma(H',Y')$,
and $\sigma(H,Y')=\sigma(H'',Y'')$ otherwise; 
\item 
if  $y'\notin Y'$ and $\sigma(H,Y')=1$,
then $\pi_i(H,Y')=\pi_i(H',Y')$ and 
$\tau_i(H,Y')=\tau_i(H',Y')$ for all $i=1,2,\cdots,n$,
and  
$D(H,Y',x')=D(H',Y',x')\cup \{y'\}$ 
and $D(H,Y',x)=D(H',Y',x)$ for all $x\in X-\{x'\}$.
\item if $y'\in Y'$ and $\sigma(H,Y')=1$,  
then 
$\pi_1=\pi_1(H,Y')$ and $\tau_1=\tau_1(H,Y')$
such that 
$y_{\tau_1}=y'$ and $x_{\pi_1}=x'$,
and $\pi_i(H,Y')=\pi_{i-1}(H'',Y'')$ and 
$\tau_i(H,Y')=\tau_{i-1}(H'',Y'')$ for all $i=2,3,\cdots,n$,
and 
$D(H,Y',x')=\emptyset$
and $D(H,Y',x)=D(H'',Y'',x)$ for all $x\in X-\{x'\}$.
\end{enumerate}
\end{lem}

Lemma~\ref{le2-2} implies that  when $\sigma(H,Y')=1$,
the outputs $\pi_i$ and $\tau_i$ are 
independent of the vertices in $Y-Y'$,
but each set $D(x_{\pi_i})$ is a subset of $Y-Y'$. 
Now we are going to show that 
when $\sigma(H,Y')=1$,
the outputs of running Algorithm A
can be determined by the following result.

\begin{pro}\relabel{prop3-1}
Let $Y'\subseteq Y$ with $\sigma(H,Y')=1$.
Then 
$\pi_i, \tau_i$ and $D(x_{\pi_i})$ 
for $i=1,2,\cdots,n$ 
can be determined by the following statements:   
\begin{enumerate}
\item for $i=1,2,\cdots,n$, $x_{\pi_i}y_{\tau_i}\in E(H)$ 
and 
$y_{\tau_i}$ is the member in $Y'\cap L(H_i)$ 
with the minimum weight $w(y_{\tau_i})$,
where $H_{i}$ denotes the subgraph of $H$ 
induced by $\bigcup_{i\le s\le n}N[x_{\pi_s}]$;
\item for $i=1,2,\cdots,n$, 
$D(x_{\pi_i})$ is the set of those 
$y\in (Y-Y')\cap N(x_{\pi_i})\cap L(H_s)$ 
such that 
$w(y)<w(y_{\tau_s})$ holds
for some $s$ with $s\le i$. 
\end{enumerate}
\end{pro}

\proof (i). 
By Lemma~\ref{le2-2} (i), 
$\pi_i$ and $\tau_i$ for $i=1,2,\cdots,n$ 
are determined by running Algorithm A with 
input $(H[X\cup Y'],Y')$.

It can  be proved by induction on $|X|$. 
The result is obvious when $|X|=1$.

Now assume that $|X|\ge 2$.
By Lemma~\ref{le2-2}, 
$\tau_1$ is determined by the fact that 
$y_{\tau_1}$ is the member in $Y'\cap L(H_1)$
with the minimum weight $w(y_{\tau_1})$
and $\pi_1$ is determined by the fact that 
$x_{\pi_1}$ is the only member in $N_{H_1}(y_{\tau_1})$.
By the inductive hypothesis, $\pi_i$ and $\tau_i$ for $i=2,\cdots,n$ 
are determined by running Algorithm A 
with the input 
$(H',Y'-\{y_{\tau_1}\})$,
where $H'=H[(X\cup Y')-\{x_{\pi_1}, y_{\tau_1}\}]$.
Thus (i) holds.
 
(ii) By Lemma~\ref{le2-2}, $\bigcup_{1\le i\le n}D(x_{\pi_i})$ 
consists of those $y\in (Y-Y')\cap L(H_s)$ 
with $w(y)<w(y_{\tau_s})$ for 
some $s: 1\le s\le n$.
Furthermore, if $y\in (Y-Y')\cap L(H_s)$ 
with $w(y)<w(y_{\tau_s})$ for 
some $s: 1\le s\le n$,
then $y\in D(x_{\pi_i})$, 
where $x_{\pi_i}$ is the only member in $N_{H_s}(y)$.
Clearly $i\ge s$ and $x_{\pi_i}$ is the only vertex 
in the set $\{x_{\pi_j}: s\le j\le n\}$ 
which is adjacent to $y$.

Hence (ii) holds.
\proofend

By Theorem~\ref{gol-theo} and Proposition~\ref{prop3-1}, 
we have the following corollary.  

\begin{cor}\relabel{cor3-0}
Let $Y'\subseteq Y$ with $\sigma(H,Y')=1$. 
Then 
\begin{enumerate}
\item $\{x_{\pi_i}y_{\tau_i}:i=1,2,\cdots,n\}$ is a 
member in $\UMX(H)$;
\item  $x_{\pi_j}y_{\tau_i}\notin E(H)$ for all $j$ with $j>i$;
\item 
if $y_{\tau_i}, y_{\tau_j}\in L(H_r)$, where $r\le \min\{i, j\}$,
then $w(y_{\tau_i})<w(y_{\tau_j})$ \iff $i<j$,
where $H_{r}$ is the subgraph of $H$ 
induced by $\bigcup_{r\le s\le n}N[x_{\pi_s}]$.
\end{enumerate} 
\end{cor}

When $\sigma(H,Y')=1$, let $M_{H,Y'}$ denote the 
subset $\{x_{\pi_i}y_{\tau_i}:i=1,2,\cdots,n\}$ of $E(H)$.
It will be shown in Corollary~\ref{cor3-2} 
that for any $T\in \sett(G)$, if $Y'=E(T)$ and $H$ is the graph $H_{G,x_0}$, then $M_{H,Y'}=M_T$.

By Corollary~\ref{cor3-0} (i),
$M_{H,Y'}\in \UMX(H)$.
Thus $\sigma(H,Y')=1$ implies that  $V(M)\cap Y\subseteq Y'$
holds for some $M\in \UMX(H)$.
Now we show that its converse statement also holds.

\begin{pro}\relabel{prop3-0}
Assume that $Y'\subseteq Y$.
Then $\sigma(H,Y')=1$ \iff 
$V(M)\cap Y\subseteq Y'$ holds for some $M\in \UMX(H)$.
\end{pro}

\proof By Corollary~\ref{cor3-0} (i),
the necessity holds.
It suffices to prove the sufficiency.

When $|X|=|Y|=1$, it is clear that 
the sufficiency holds.
Assume that the sufficiency holds when $2\le |X|+|Y|<r$.
Now consider the case that $|X|+|Y|=r$
and 
assume that there exists $M\in \UMX(H)$
with $V(M)\cap Y\subseteq Y'$. 

As $\UMX(H)\ne \emptyset$, 
by Theorem~\ref{gol-theo}, 
$L(H)\cap Y\ne \emptyset$. 
Let $y'$ be the member in 
$L(H)\cap Y$ such that $w(y')$ is the minimum.
If $y'\notin Y'$, then $M\in \UMX(H')$
with $V(M)\cap (Y-\{y'\})\subseteq Y'$, 
where $H'=H-y'$, 
and by the inductive hypothesis,  $\sigma(H',Y')=1$ holds.
If $y'\in Y'$, then $M-\{x'y'\}\in \UMX(H'')$,
where $H''=H-\{x',y'\}$
and $x'$ is the only member in $N_H(y')$,  
and by the inductive hypothesis,  $\sigma(H'',Y'')=1$ holds,
where $Y''=Y'-\{y'\}$.
In both cases, Lemma~\ref{le2-2} implies that 
 $\sigma(H,Y'')=1$.

Hence the sufficiency holds. 
\proofend


\subsection{A member of $\BPX(H)$ 
when $\sigma(H,Y')=1$
\relabel{sect3-3}}

When $\sigma(H,Y')=1$,
a special member of $\BPX(H)$ can be determined 
by the sets $D(H,Y',x)$'s.

\begin{pro}\relabel{prop2-44}
For any $Y'\subseteq Y$ with $\sigma(H,Y')=1$,
the function $f:X\rightarrow \N_0$ 
determined by $f(x)=|D(H,Y',x)|$ for all 
$x\in X$ is a member in $\BPX(H)$.
\end{pro}

\proof 
We prove it by induction on $|X|+|Y|$.
The result is obvious when $|X|=|Y|=1$ by 
Proposition~\ref{prop2-11}. 
Assume that the result holds when $2\le |X|+|Y|<r$.
Now consider the case that $|X|+|Y|=r$.

As $\sigma(H,Y')=1$,
$L(H)\cap Y\ne \emptyset$. 
Let $y'$ be the member in 
$L(H)\cap Y$ such that $w(y')$ is the minimum.
Let $x'$ be the only member in $N_H(y')$.

First consider the case that $y'\notin Y'$.
By the inductive hypothesis,
the function $g:X\rightarrow \N_0$ 
defined by $g(x)=|D(H',Y',x)|$ for all 
$x\in X$ is a member in $\BPX(H')$, 
where $H'=H-y'$.
By Corollary~\ref{cor2-02}(ii),
the function $f:X\rightarrow \N_0$ 
defined by $f(x')=g(x')+1$
and $f(x)=g(x)$ for all $x\in X-\{x'\}$ 
is a member in $\BPX(H)$.
By Lemma~\ref{le2-2}(i),
$f(x)=|D(H,Y',x)|$ for all $x\in X$.
Thus the result holds in this case.

Now consider the case that $y'\in Y'$.
Then $\sigma(H'',Y'')=\sigma(H,Y')=1$ 
by Lemma~\ref{le2-2} (ii),
where $Y''=Y'-\{y'\}$
and $H''=H-\{x',y'\}$.
By the inductive hypothesis, the function $g:X-\{x'\}\rightarrow \N_0$
defined by 
$g(x)=|D(H'',Y'',x)|$ for all $x\in X-\{x'\}$ 
is a member in $\BP_{X'}(H'')$, where $X'=X-\{x'\}$.
By Corollary~\ref{cor2-02}(iii),
the function $f:X\rightarrow \N_0$ 
defined by $f(x')=0$
and $f(x)=g(x)$ for all $x\in X-\{x'\}$ 
is a member in $\BPX(H)$.
By Lemma~\ref{le2-2}(ii),
$f(x)=|D(H,Y',x)|$ for all $x\in X$.
Thus the result also holds in this case.

Hence the result holds.
\proofend

\subsection{Outputs of running Algorithm A for 
$H_{G,x_0}$\relabel{sect3-4}}

In the next two subsections, 
we will consider the special bipartite graph 
$H_{G,x_0}$. 

Note that the weight function $w:E\rightarrow \N_0$ for
edges of $G$ 
is also the weight function for members of $Y$ in 
$H_{G,x_0}$, which is 
used in running Algorithm $A$ with input $(H_{G,x_0},Y')$,
where $Y'\subseteq Y=E$.
If $\sigma(H_{G,x_0},Y')=1$, simply write $\pi_i=\pi_i(H_{G,x_0},Y')$, 
$\tau_i=\tau_i(H_{G,x_0},Y')$ and $D(x_{\pi})=D(H_{G,x_0},Y',x_{\pi})$
for $i=1,2,\cdots,n$.

The next result follows from
Propositions~\ref{prop3-1} and~\ref{prop3-0}.

\begin{pro}\relabel{prop3-5}
$\sigma(H_{G,x_0},Y')=1$ \iff $G[Y']$ 
is a connected and spanning subgraph of $G$.
Furthermore, if  $\sigma(H_{G,x_0},Y')=1$, then,  
for $i=1,2,\cdots,n$, 
\begin{enumerate}
\item  
$y_{\tau_i}$ is the edge in $Y'\cap E_G(V_i, V-V_i)$ 
with $w(y_{\tau_i})\le w(y')$ for all 
$y'\in Y'\cap E_G(V_i, V-V_i)$
and $x_{\pi_i}$ is the vertex in $V-V_i$ incident with 
$y_{\tau_i}$, where $V_i=\{x_0\}\cup \{x_{\pi_s}:1\le s<i\}$;

\item 
$D(x_{\pi_i})$ is the set of those edges $y\in Y-Y'$ incident with $x_{\pi_i}$
such that $y\in E_G(V_s, V-V_s)$ and 
$w(y)<w(y_{\tau_s})$ hold for some $s\le i$.

\end{enumerate}
\end{pro}

By Lemma~\ref{le2-3} and Proposition~\ref{prop3-5} (i), 
for any $T\in \sett(G)$, 
we have the following relation on $M_T$ and $M_{H_{G,x_0},Y'}$, where $Y'=E(T)$. 

\begin{cor}\relabel{cor3-2}
For any $T\in \sett(G)$, 
$M_T=M_{H_{G,x_0},Y'}=\{x_{\pi_i}y_{\tau_i}:i=1,2,\cdots,n\}$,
where $Y'=E(T)$.
\end{cor}

For example, let $G=(V,E)$ be the graph shown 
in Figure~\ref{S3-f1} (a) and 
$Y'$ be a subset of $E$ with $G[Y']$ shown in Figure~\ref{S3-f1} (b),
where each number beside an edge $e$ is its weight $w(e)$.
As $G[Y']$ is a spanning tree of $G$,
Proposition~\ref{prop3-5} implies that 
$\sigma(H_{G,x_0},Y')=1$.

\begin{figure}[htbp]
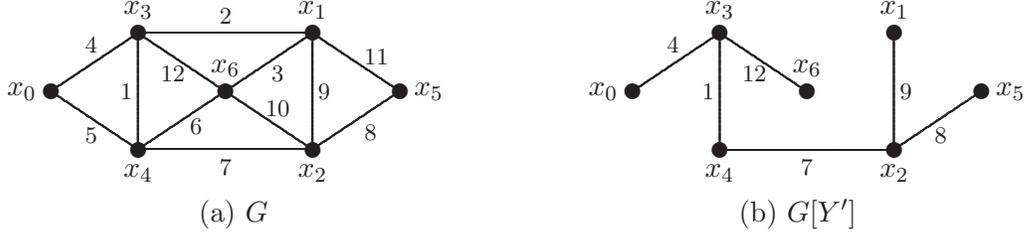

\centering

\beginpicture 
\setcoordinatesystem units <2.2pt,2.2pt> 
\setplotarea x from -135 to 45, y from 0 to 20 
\plotsymbolspacing=.2pt

\put {\cirs} at -120 10  \put {$x_0$} at  -125 10
\put {\cirs} at -90 10  \put {$x_6$} at  -90 14
\put {\cirs} at -60 10  \put {$x_5$} at  -55 10
\put {\cirs} at -105  0  \put {$x_4$} at  -105 -4
\put {\cirs} at -75 0  \put {$x_2$} at  -75 -4
\put {\cirs} at -105  20  \put {$x_3$} at  -105 24
\put {\cirs} at -75 20  \put {$x_1$} at  -75 24

\plot -120 10 -105 0 -75 0  -60 10 -75 20 -105 20  -120 10 / 

\plot -75 20 -75 0 / 
\plot -105 20 -105 0 / 

\plot -90 10 -75 20 /
\plot -90 10 -75 0 /
\plot -90 10 -105 20 /
\plot -90 10 -105 0 /

\put{\footnotesize $1$} at -107 10 
\put{\footnotesize $4$} at -113 18
\put{\footnotesize $5$} at -113 2.5  
\put{\footnotesize $9$} at -73 10
\put{\footnotesize $11$} at -64 16
\put{\footnotesize $10$} at -81 7
\put{\footnotesize $3$} at -81 13
\put{\footnotesize $6$} at -95 4
\put{\footnotesize $12$} at -99 13
\put{\footnotesize $8$} at -65 3
\put{\footnotesize $2$} at -90 23 
\put{\footnotesize $7$} at -90 -3



\put {\cirs} at -20 10  \put {$x_0$} at  -25 10
\put {\cirs} at 10 10  \put {$x_6$} at  10 14
\put {\cirs} at 40 10  \put {$x_5$} at  45 10
\put {\cirs} at -5  0  \put {$x_4$} at  -5 -4
\put {\cirs} at 25 0  \put {$x_2$} at  25 -4
\put {\cirs} at -5  20  \put {$x_3$} at  -5 24
\put {\cirs} at 25 20  \put {$x_1$} at  25 24

\plot -5 20 -20 10 / 
\plot 25 20 25 0 / 
\plot -5 20 -5 0 / 
\plot 10 10 -5 20 /

\plot 25 0 -5 0 /
\plot  40 10 25 0 /

\put{\footnotesize $1$} at -7 10 
\put{\footnotesize $4$} at -13 18
\put{\footnotesize $9$} at 27 10
\put{\footnotesize $12$} at 1 13
\put{\footnotesize $8$} at 33 2.5
\put{\footnotesize $7$} at 10 -3 

\endpicture 

\vspace{0.2 cm}

(a) $G$ \hspace{6 cm} (b) $G[Y']$

\caption{$G$ and a spanning tree $G[Y']$}
\relabel{S3-f1}
\end{figure}

By Proposition~\ref{prop3-5} (i), 
$y_{\tau_1}, \cdots, y_{\tau_6}$ are the following edges respectively:
$$
x_0x_3, x_3x_4, x_4x_2, x_2x_5,x_2x_1, x_3x_6, 
$$
and $x_{\pi_1}, \cdots, x_{\pi_6}$ 
are the vertices $x_3, x_4, x_2, x_5,x_1, x_6$
respectively.
By applying Proposition~\ref{prop3-5} (ii), 
we have 
$$
D(x_3)=D(x_4)=D(x_2)= D(x_5)=\emptyset,
D(x_1)=\{x_3x_1\},
D(x_6)=\{x_4x_6,x_2x_6,x_1x_6\}.
$$

The next result considers the special case 
that $Y'=E(T)$ for a given $T\in \sett(G)$.
It will be applied for 
proving Theorem~\ref{theo5-1}.

Let $P_{i,j}$ denote the unique path in $T$ 
connecting vertices $x_{\pi_i}$ and $x_{\pi_j}$.

\begin{pro}\relabel{prop3-6}
Let  $T\in \sett(G)$ 
and $Y'=E(T)$. 
Then 
\begin{enumerate}
\item for  $i=1,2,\cdots,n$,
$G[E_i]$ is a tree with 
vertex set $\{x_{\pi_s}:0\le s\le i\}$,
where $E_i=\{y_{\tau_s}: 1\le s\le i\}$
and $\pi_0=0$; 

\item for $i=1,2,\cdots,n$,  
$y_{\tau_i}$ is incident with $x_{\pi_i}$
and is an edge on the path $P_{0,i}$;
 
\item 
if $x_{\pi_i}$ is a vertex on 
the path $P_{0,j}$, then $i\le j$ holds;

\item for any integers $1\le i,j\le n$, 
if $\max\{b(y_{\tau_i}), b(y_{\tau_j})\}<\min\{i,j\}$,
then 
$w(y_{\tau_i})<w(y_{\tau_j})$ \iff $i<j$,
where $b(y_{\tau_j})$ is the number $s$ 
such that $x_{\pi_s}$ is the end of $y_{\tau_j}$ in $G$
different from $x_{\pi_j}$.
\end{enumerate}
\end{pro}

\proof (i) follows from Proposition~\ref{prop3-5} (i).

(ii) and (iii) follow directly from result (i).

(iv). Let $r=\max\{b(y_{\tau_i}), b(y_{\tau_j})\}$.
As $r<\min\{i,j\}$, 
both $y_{\tau_i}$ and $y_{\tau_j}$ are 
members in the set $Y'\cap E_G(V_k,V-V_k)$ 
for all $k$ with $r<k\le \min\{i,j\}$,
where $V_k=\{x_{\pi_t}:k\le t\le n\}$.
By Proposition~\ref{prop3-5} (i),
$w(y_{\tau_i})<w(y_{\tau_j})$ \iff 
$y_{\tau_i}$ is selected before 
$y_{\tau_j}$, i.e., $i<j$.
Thus (iv) holds.
\proofend

\subsection{The minimum spanning tree \relabel{sect3-5}}

The {\it minimum spanning tree} of $G$
with respect to $w$ 
is the spanning tree $T_0$ of $G$ such that 
$w(T_0)<w(T)$ holds for all $T\in \sett(G)-\{T_0\}$,
where $w(T)=\sum_{e\in E(T)}w(e)$. 
In this subsection, we  show that 
the minimum spanning tree of $G$ 
is determined by the outputs $y_{\tau_i}$'s 
of running Algorithm A with input $(H_{G,x_0},E(G))$.
But this property cannot be extended to 
all bipartite graphs. 

Prim's algorithm \cite{prim} is 
a well-known algorithm 
of determining the minimum spanning tree 
of a connected multigraph.
The way of choosing edges 
of the minimum spanning tree in $G$ 
by Prim's algorithm (see \cite{bru,west})
is actually the same as the way of 
determining edges 
$y_{\tau_1},\cdots,y_{\tau_n}$ by Proposition~\ref{prop3-5} (i).
Thus the next result follows from 
Proposition~\ref{prop3-5} (i) and Prim's algorithm.

\begin{cor}\relabel{cor3-02}
For any $Y'\subseteq E$, 
if $G[Y']$ is connected and spanning, 
then $E(T_0)=\{y_{\tau_1},\cdots,y_{\tau_n}\}$ 
for the minimum spanning tree $T_0$ of $G[Y']$.
\end{cor}

For any $Y''\subseteq E$ with $Y'\subset Y''$, 
when do $G[Y']$ and $G[Y'']$ have the same 
the minimum spanning tree?

\begin{theo}\relabel{theo3-2}
Let $T_0$ be the minimum spanning tree of $G[Y']$.
For any $Y''\subseteq E$ with $Y'\subseteq Y''$, 
$T_0$ is the minimum spanning tree of $G[Y'']$ \iff 
$\left (\bigcup_{1\le i\le n}D(x_{\pi_i})\right )\cap Y''=\emptyset$.
\end{theo}

\proof
It suffices to show that the two statements below hold: 

(a) if $\left (\bigcup_{1\le i\le n}D(x_{\pi_i})\right )\cap Y''=\emptyset$,
then $T_0$ is the minimum spanning tree of $G[Y'']$;

(b) if $\left (\bigcup_{1\le i\le n}D(x_{\pi_i})\right )\cap Y''\ne \emptyset$,
then $T_0$ is not the minimum spanning tree of $G[Y'']$.

Assume that $\left (\bigcup_{1\le i\le n}D(x_{\pi_i})\right )\cap Y''=\emptyset$.
By Proposition~\ref{prop3-5} (ii), 
$\bigcup_{1\le i\le n}D(x_{\pi_i})$ is the set of those 
edges $y\in Y-Y'$ such that 
$y\in E_G(V_s,V-V_s)$ and $w(y)\le w(y_{\tau_s})$ 
hold for some $s$ with $1\le s\le n$, 
where $V_s=\{x_{\pi_t}: s\le t\le n\}$.
As $\left (\bigcup_{1\le i\le n}D(x_{\pi_i})\right )\cap Y''=\emptyset$,
by Proposition~\ref{prop3-5} (i), 
$y_{\tau_i}$ is the edge in  $E_{G[Y'']}(V_i,V-V_i)$ 
such that 
$w(y_{\tau_i})<w(y)$ holds for all edges 
$y\in E_{G[Y'']}(V_i,V-V_i)-\{y_{\tau_i}\}$
for each $i=1,2,\cdots,n$.
By Prim's algorithm, 
$E(T_0)=\{y_{\tau_i}:i=1,2,\cdots,n\}$ is the edge set 
of the minimum spanning tree of $G[Y'']$.
Hence (a) holds.

Now consider the case that 
$\left (\bigcup_{1\le i\le n}D(x_{\pi_i})\right )
\cap Y''\ne \emptyset$.
By Corollary~\ref{cor3-02}, 
the edge set of the minimum spanning tree $T_0$ of 
$G[Y']$ is 
$\{y_{\tau_1}, y_{\tau_2},\cdots, y_{\tau_n}\}$.
By Prim's algorithm, 
the edges of $T_0$ can be chosen 
in the order $y_{\tau_1}, y_{\tau_2},\cdots, y_{\tau_n}$.
By the assumption, 
there exists 
$y_0\in \left (\bigcup_{1\le i\le n}D(x_{\pi_i})\right )
\cap Y''$.
By Proposition~\ref{prop3-5} (ii),
$y_0\in E_{G[Y'']}(V_s, V-V_s)$ and 
$w(y_0)<w(y_{\tau_s})$ hold for some $s$ with $1\le s\le n$.
By Prim's algorithm again, 
$y_0$ is chosen  
as an edge of the minimum spanning tree of $G[Y'']$
at the step after all edges in $\{y_{\tau_t}: 1\le t<s\}$ 
are selected, implying that 
$T_0$ is not the minimum spanning tree of 
$G[Y'']$.

Hence (b) also holds. 
\proofend

For any $M\in \UMX(H)$, let $w(M)=\sum_{y\in V(M)\cap Y} w(y)$.
A member $M_0$ in $\UMX(H)$ is called a
{\it minimum member} in $\UMX(H)$ 
if $w(M_0)\le w(M)$ holds for all 
$M\in \UMX(H)$. 
By Corollary~\ref{cor3-02},  
$\{x_{\pi_i}y_{\tau_i}:i=1,2,\cdots,n\}$ is 
the unique minimum member of $\UMX(H_{G,x_0})$.
However, this result does not hold all bipartite graphs $H$. 
An example is shown in Figure~\ref{S3-f11}.

Let $H_0$ be the bipartite graph shown
in Figure~\ref{S3-f11},
where any vertex with an order pair $(y_i,w_i)$ beside 
is vertex $y_i$ with $w(y_i)=w_i$.
Running Algorithm A with input $(H_0,Y_0)$,
where $Y_0=\{y_1,y_2,y_3,y_4\}$, outputs $\pi_i=\tau_i=i$ for $i=1,2,3$. 
But $M_0=\{x_iy_i:i=1,2,3\}$ is 
not the minimum member  of $\UMX(H_0)$, 
as $M_1=\{x_2y_2,x_3y_3,x_1y_4\}\in \UMX(H_0)$
and 
$$
w(M_1)=w(y_2)+w(y_3)+w(y_4)<w(y_1)+w(y_2)+w(y_3)=w(M_0).
$$

\begin{figure}[htbp]
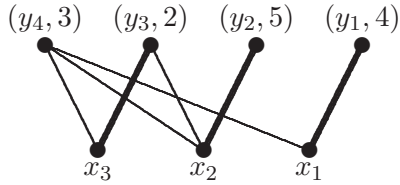

\centering

\begin{center}
\beginpicture 
\setcoordinatesystem units <2pt,2pt> 
\setplotarea x from -105 to 100, y from -4 to 20 
\plotsymbolspacing=.18pt 

\put {\cirs} at 20 0  \put {$x_1$} at  20 -4
\put {\cirs} at 0 0  \put {$x_2$} at  0 -4
\put {\cirs} at -20 0  \put {$x_3$} at  -20 -4

\put {\cirs} at 30 20  \put {$(y_1,4)$} at  30 25
\put {\cirs} at 10 20  \put {$(y_2,5)$} at  10 25
\put {\cirs} at -10 20  \put {$(y_3,2)$} at  -10 25
\put {\cirs} at -30 20  \put {$(y_4,3)$} at  -30 25

\plot 30 20 20 0 -30 20 0 0 10 20 / 
\plot 0 0 -10 20 -20 0 -30 20 / 

\plot 30.4 20 20.4 0 /  \plot 29.6 20 19.6 0 / 

\plot 10.4 20 0.4 0 /  \plot 9.6 20 -0.4 0 / 

\plot -10.4 20 -20.4 0 /  \plot -9.6 20 -19.6 0 / 

\endpicture 
\end{center}

\vspace{-0.3 cm}

\caption{A bipartite graph $H_0$}
\relabel{S3-f11}
\end{figure}

\begin{prob}\relabel{prob3-1}
For any bipartite graph $H$ with a bipartition $(X,Y)$
and $\UMX(H)\ne \emptyset$,
determine the minimum member of $\UMX(H)$.
\end{prob}

\resection{Bijection $\psi_H$ from $\UMX(H)$ to $\BPX(H)$
\relabel{sect04}}

For any $M\in \UMX(H)$,
let $\psi_H(M)=f$, where $f$ is the mapping $f:X\rightarrow \N_0$
defined by 
$f(x_i)=|D(H,Y\cap V(M),x_i)|$ for each $x_i\in X$.
By Propositions~\ref{prop3-0} and~\ref{prop2-44}, 
$\psi_H$ is a mapping from $\UMX(H)$ to $\BPX(H)$.
By its definition, 
an interpretation of $\psi_H$ is given by 
Proposition~\ref{prop3-1}.
We are now going to show that $\psi_H$ is a bijection.

\begin{theo}\relabel{B-bi}
The mapping $\psi_H:\UMX(H)\rightarrow \BPX(H)$ defined 
above is a bijection from $\UMX(H)$ to $\BPX(H)$. 
\end{theo}

\proof We first prove that $\psi_H$ is injective by induction
on $|X|+|Y|$. 
When $|X|=|Y|=1$, the conclusion is obvious,
as $\UMX(H)$ has at most one member. 
Assume that it holds when $|X|+|Y|<k$, where $k\ge 3$.
Now consider the case that $|X|+|Y|=k$.

Assume that $\UMX(H)\ne \emptyset$. 
By Theorem~\ref{gol-theo}, 
$Y\cap L(H)\ne \emptyset$.
Assume that $y'$ is the member 
in $Y\cap L(H)$ such that $w(y')$ is the minimum. 
Let $x'$ be the only member in $N_H(y')$. 

Let $M_1$ and $M_2$ be distinct members in $\UMX(H)$
and $Y_i=V(M_i)\cap Y$ for $i=1,2$.
If $Y_1=Y_2$, then $V(M_1)=V(M_2)$,
implying that $M_1=M_2$ by the definition of 
UR-matchings. Thus $Y_1\ne Y_2$.
Let $f_i(x)=|D(H,Y_i,x)|$ for $i=1,2$ and all $x\in X$.
We shall show that $f_1\ne f_2$ in the three cases below.

\noindent {\bf Case 1}: $y'\in Y_1-Y_2$ or $y'\in Y_2-Y_1$.

Assume that $y'\in Y_1-Y_2$.
By Lemma~\ref{le2-2}, $D(H,Y_1,x')=\emptyset$ 
while $y'\in D(H,Y_2,x')$. 
Thus $f_1(x')<f_2(x')$. 

\noindent {\bf Case 2}: $y'\notin Y_1\cup Y_2$.
 

In this case, $M_i\in \UMX(H')$ for $i=1,2$, where $H'=H-y'$. 
By the inductive hypothesis, $\psi_{H'}$ is an injective mapping from $\UMX(H')$
to $\BPX(H')$, implying that  
$|D(H',Y_1,x)|\ne |D(H',Y_2,x)|$ for some $x\in X$.
By Lemma~\ref{le2-2}(i), for each $i=1,2$, 
$D(H,Y_i,x')=D(H',Y_i,x')\cup \{y'\}$ 
and $D(H,Y_i,x)=D(H',Y_i,x)$  for all $x\in X-\{x'\}$,
implying that 
$|D(H,Y_1,x)|\ne |D(H,Y_2,x)|$ for some $x\in X$,
i.e., $f_1\ne f_2$.

\noindent {\bf Case 3}: $y'\in Y_1\cap Y_2$.

By Lemma~\ref{le2-2}(ii), for $i=1,2$,
$D(H,Y_i,x')=\emptyset$  
and $D(H,Y_i,x)=D(H'',Y_i',x)$  for all 
$x\in X'=X-\{x'\}$, where $H''=H-\{x',y'\}$ 
and $Y_i'=Y_i-\{y'\}$. 
Note that $Y'_i=Y\cap V(M'_i)$ for $i=1,2$,
where $M'_i=M_i-\{x'y'\}$.
As $M_1\ne M_2$, we have $M'_1\ne M'_2$.
By the inductive hypothesis, 
$|D(H'',Y_1',x)|\ne |D(H'',Y_2',x)|$
for some $x\in X'$,
implying that 
$|D(H,Y_1,x)|\ne |D(H,Y_2,x)|$.
Thus, $f_1\ne f_2$ in this case.

Therefore $\psi_H$ is injective.

It remains to prove that $\psi_H$ is surjective,
i.e., the following statement ``for any $f\in \BPX(H)$, 
there exists $M\in \UMX(H)$ with $\psi_H(M)=f$" holds.
We prove this statement by induction on the value of $|X|+|Y|+\sum_{x\in X}f(x)$,
where $f\in \BPX(H)$. 
Observe that $|X|+|Y|+\sum_{x\in X}f(x)\ge 2$.
When $|X|+|Y|+\sum_{x\in X}f(x)=2$, 
we have $|X|=|Y|=1$ and $f(x)=0$ for the only member $x\in X$,
implying that $H\cong K_2$ and $\psi_H(M)=f$ holds,
where $M=E(H)$. 

Assume that the above statement holds 
for any bipartite graph $H'$ with 
a bipartition $(X',Y')$ 
and any $f'\in \BP_{X'}(H')$ such that 
$|X'|+|Y'|+\sum_{x\in X'}f'(x)<r$, where $r\ge 3$.
Now we suppose that $H$ is a bipartite graph with 
a bipartition $(X,Y)$ and $f\in \BPX(H)$ 
such that $|X|+|Y|+\sum_{x\in X}f(x)=r$.

As $\BPX(H)\ne \emptyset$, 
by Proposition~\ref{prop2-12},
we have $\UMX(H)\ne \emptyset$
and 
$Y\cap L(H)\ne \emptyset$.
Assume that $y'$ is the member 
in $Y\cap L(H)$ such that $w(y')$ is the minimum
and $x'$ is the only member in $N_H(y')$.
We shall prove in the two cases below that 
$\psi_H(M)=f$ holds for some $M\in \UMX(H)$.

\noindent {\bf Case 1'}: $f(x')=0$.

Let $H''=H-\{x',y'\}$ and $g=f|_{X'}$, where $X'=X-\{x'\}$.
By Corollary~\ref{cor2-02}(iii), $g\in \BP_{X'}(H'')$.
By the inductive hypothesis, there exists $M'\in \UM_{X'}(H'')$ 
such that $\psi_{H''}(M')=g$,
i.e., $g(x)=|D(H'',V(M')\cap Y,x)|$ for all $x\in X'$.
It is clear that
$M=M'\cup \{x'y'\}\in \UMX(H)$.
By Lemma~\ref{le2-2}(ii), 
$D(H,Y',x')=\emptyset$ and 
$D(H,Y',x)=D(H'',Y'',x)$ for all $x\in X-\{x'\}$,
where $Y''=V(M')\cap Y$ and $Y'=Y''\cup \{y'\}=V(M)\cap Y$.
Thus $f(x')=0=|D(H,Y',x')|$ and 
$f(x)=g(x)=|D(H'',Y'',x)|=|D(H,Y',x)|$
for all $x\in X-\{x'\}$, 
implying that $\psi_H(M)=f$.

\noindent {\bf Case 2'}: $f(x')>0$.

Let $H'=H-\{y'\}$ and $g=f_{(x'\downarrow 1)}$.
By Corollary~\ref{cor2-02}(ii), $g\in \BPX(H')$.
By the inductive hypothesis, there exists $M\in \UMX(H')$ 
such that $\psi_{H'}(M)=g$,
i.e., $g(x)=|D(H',Y',x)|$ for all $x\in X$,
where $Y'=V(M)\cap Y$.
By Lemma~\ref{le2-2}(i), 
$D(H,Y',x')=D(H',Y',x')\cup \{y'\}$ and 
$D(H,Y',x)=D(H',Y',x)$ for all $x\in X-\{x'\}$.
Thus $f(x')=g(x')+1=|D(H',Y',x')|+1
=|D(H,Y',x')|$ and 
$f(x)=g(x)=|D(H',Y',x)|=|D(H,Y',x)|$
for all $x\in X-\{x'\}$, 
implying that $\psi_H(M)=f$.
\proofend

For any $T\in \sett(G)$, define  $\phi_G(T)=\psi_{H_{G,x_0}}(M_T)$.
By Theorem~\ref{B-bi}, Corollary~\ref{cor3-2}
and Proposition~\ref{prop2-4},
$\phi_G$ is a bijection from $\sett(G)$ 
to $\GP(G,x_0)$.
By Proposition~\ref{prop3-5}, $\phi_G$ 
can be interpreted  by the following result,
which first appeared in~\cite{che}.

\begin{cor}\relabel{cor3-3}
Let $T\in \sett(G)$.
Assume that vertices $x_{\pi_1},  x_{\pi_2},\cdots, x_{\pi_n}$ 
and egdes $y_{\tau_1},  y_{\tau_2},\cdots, y_{\tau_n}$ 
of $G$ 
are determined by Proposition~\ref{prop3-5} (i),
where $Y'=E(T)$.  
If $f=\phi_G(T)$, then, for $i=1,2,\cdots,n$, 
$f(x_{\pi_i})$ is the number of those edges 
$y'\in E(G)-E(T)$ incident with 
$x_{\pi_i}$ and some $x_{\pi_j}$,
where $0\le j<i$,
with $w(y')<\max_{j<s\le i} w(y_{\tau_s})$.
\end{cor}

For example, if $G$ is the graph shown in Figure~\ref{S3-f1} (a)
and $T$ is the spanning tree in Figure~\ref{S3-f1} (b), 
then $\phi_G(T)$ is the mapping $f\in \GP(G,x_0)$ given below: 
$$
f(x_2)=f(x_3)=f(x_4)=f(x_5)=0, f(x_1)=1, f(x_6)=3.
$$

\resection{Interpret B-parking functions}
\relabel{sect5}

Theorem~\ref{B-bi} shows that 
the mapping $\psi_H:\UMX(H)\rightarrow \BPX(H)$ 
defined by $\psi_H(M)=f$ is a bijection, where 
$f(x)=|D(H,V(M)\cap Y,x)|$ for all $x\in X$.
In this section, 
assume that  $M\in \UMX(H)$ and $Y'=V(M)\cap Y$, unless otherwise stated.
Also assume that $\pi_i=\pi_i(H,Y')$, 
$\tau_i=\tau_i(H,Y')$ and $D(x_{\pi_i})=D(H,Y',x_{\pi_i})$. 
In this section, we will give an interpretation 
for $f$ different from Proposition~\ref{prop3-1} (ii).

In Subsection~\ref{sect5-1}, we define a unique path 
$P_{(H,M)}(y)$ in $H$ for each $y\in Y-Y'$ 
with respect to $M$. 
In Subsection~\ref{sect5-2}, we 
introduce the concept ``externally B-active members
with respect to $M$ in $H$" by comparing 
$w(y)$ with $w(y')$ for all those $y'\in Y$ 
which are in the path $P_{(H,M)}(y)$.
In Subsection~\ref{sect5-3}, 
we show that 
$\bigcup_{x\in X}D(H,Y',x)$ 
is exactly the set of those members in $Y-Y'$ which are not 
externally B-active with respect to $M$ in $H$.
In particular, 
$D(H,Y',x_{\pi_i})$ is the set of 
those members $y$ in $\left ((Y-Y')\cap N_H(x_{\pi_i})\right )-\bigcup_{s>i}N_H(x_{\pi_s})$ which are not externally B-active with respect to $M$ in $H$, where $Y'=V(M)\cap Y$.
Finally, in Subsection~\ref{gen-func},
we introduce a generating function 
$\Omega(H;x,y,z)$ for the members 
in $\UM(H)$ with three variables.
Particularly, 
$\Omega(H_{G,x_0};x,y,0)$ is the 
Tutte polynomial $T_G(x,y)$.

\subsection{The path $P_{(H,M)}(y)$ for each $y\in Y-Y'$}
\relabel{sect5-1}

By the definition of $\pi_i$ and $\tau_i$ for $i=1,2,\cdots,n$,
we have $Y'=\{y_{\tau_i}:i=1,2,\cdots,n\}$
and $M=M_{H,Y'}=\{x_{\pi_i}y_{\tau_i}: i=1,2,\cdots,n\}$.
For any vertex $y\in Y$ and any integer $j\ge 1$, 
let $n_j(y)=0$ if $j>d_H(y)$,
and let $n_j(y)$ be the $j$'th largest integer $s$ such that $x_{\pi_s}\in N(y)$ otherwise.
In other words, 
$n\ge n_1(y)>n_2(y)>\cdots>n_{d_H(y)}(y)>n_j(y)=0$ for all $j>d_H(y)$
and $N(y)=\{x_{\pi_s}: s\in \{n_1(y), \cdots, n_{d_H(y)}(y)\}\}$.

Clearly $n_1(y_{\tau_i})=i$ for all $i=1,2,\cdots,n$ by 
Corollary~\ref{cor3-0} (i) and (ii).
By Proposition~\ref{prop3-1}, 
$D(H,Y',x_{\pi_i})\subseteq \{y:Y-Y', n_1(y)=i\}$.

For any $y\in Y-Y'$, 
let $P_{(H,M)}(y)$ be the following maximal $M$-alternating path 
in $H$ with $y$ as one end:
$$
P_{(H,M)}(y): y x_{\pi_{j_1}}y_{\tau_{j_1}}\cdots 
x_{\pi_{j_t}}y_{\tau_{j_t}}
$$
where 
$j_1=n_1(y)$, $j_{i}=n_2(y_{\tau_{j_{i-1}}})>0$ 
for all $i=2,3,\cdots,t$
and $n_2(y)<j_t$, as shown in Figure~\ref{f7}.
Thus $j_1>j_2>\cdots>j_t>n_2(y)$.
By the maximality of $P_{(H,M)}(y)$, 
$n_2(y)\ge n_2(y_{\tau_{j_t}})\ge 0$.
Clearly that the path $P_{(H,M)}(y)$ is unique for each $y$.

\begin{figure}[htbp]
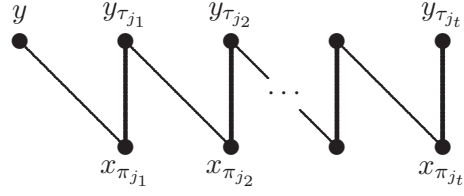

\centering

\beginpicture 
\setcoordinatesystem units <2pt,2pt> 
\setplotarea x from -100 to 100, y from 0 to 20 
\plotsymbolspacing=.2pt

\put {\cirs} at  -20 0
\put {\cirs} at  0 0
\put {\cirs} at  20 0
\put {\cirs} at  40 0

\put {\cirs} at  -40 20
\put {\cirs} at  -20 20
\put {\cirs} at  0 20
\put {\cirs} at  20 20
\put {\cirs} at  40 20

\plot -20.2 0 -20.2 20 /
\plot -19.8 0 -19.8 20 /

\plot 0.2 0 0.2 20 /
\plot -0.2 0 -0.2 20 /

\plot 20.2 0 20.2 20 /
\plot 19.8 0 19.8 20 /

\plot 40.2 0 40.2 20 /
\plot 39.8 0 39.8 20 /

\plot -20 0 -20 20 /
\plot 0 0 0 20 /
\plot 20 0 20 20 /
\plot 40 0 40 20 /
\plot -40 20 -20 0 /
\plot -20 20 0 0 /
\plot 20 20 40 0 /
\plot 0 20 7 13 /
\plot 20 0 13 7 /
\put {$\cdots$} at  10 10

\put {$x_{\pi_{j_1}}$} at  -20 -5 
\put {$x_{\pi_{j_2}}$} at  0 -5
\put {$x_{\pi_{j_t}}$} at  40 -5

\put {$y_{\tau_{j_1}}$} at  -20 25 
\put {$y_{\tau_{j_2}}$} at  0 25
\put {$y_{\tau_{j_t}}$} at  40 25
\put {$y$} at  -40 25

\endpicture

\caption{$P_{(H,M)}(y):yx_{\pi_{j_1}}y_{\tau_{j_1}}\cdots 
x_{\pi_{j_t}}y_{\tau_{j_t}}$}
\relabel{f7}
\end{figure}

For example, if $H$ is the bipartite graph 
shown in Figure~\ref{f5} with $w(y_i)=i$ for all $i$
and $M=\{x_iy_i: i=1,2,\cdots,5\}\in \UMX(H)$,
then $Y'=\{y_i:i=1,2,\cdots,5\}$, 
$\pi_i=\tau_i=i$ for $i=1,2,\cdots,5$. 
Note that $y_6$ is the only vertex in 
$Y-Y'$.
As $n_1(y_6)=5$, $n_2(y_5)=4$, $n_2(y_4)=3$
and $n_2(y_3)=2=n_2(y_6)$,  the path $P_{(H,M)}(y_6)$ is 
$$
P_{(H,M)}(y_6): y_6x_5y_5x_4y_4x_3y_3.
$$

\begin{figure}[htbp]
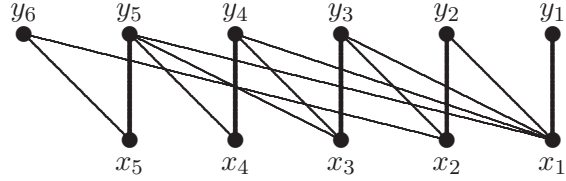

\centering

\beginpicture 
\setcoordinatesystem units <2pt,2pt> 
\setplotarea x from -110 to 100, y from 0 to 20 
\plotsymbolspacing=.2pt 

\put {\cirs} at  -60 20  \put {$y_6$} at  -60 24
\put {\cirs} at  -40 20  \put {$y_5$} at  -40 24
\put {\cirs} at  -20 20 \put {$y_4$} at  -20 24
\put {\cirs} at  0 20 \put {$y_3$} at  0 24
\put {\cirs} at  20 20 \put {$y_2$} at  20 24
\put {\cirs} at  40 20  \put {$y_1$} at  40 24

\put {\cirs} at  -40 0 \put {$x_5$} at  -40 -5
\put {\cirs} at  -20 0 \put {$x_4$} at  -20 -5
\put {\cirs} at  0 0  \put {$x_3$} at  0 -5
\put {\cirs} at  20 0 \put {$x_2$} at  20 -5
\put {\cirs} at  40 0 \put {$x_1$} at  40 -5

\plot 40 0 40 20 / \plot 39.8 0 39.8 20 / \plot 40.2 0 40.2 20 /

\plot 20 0 20 20 / \plot 19.8 0 19.8 20 / \plot 20.2 0 20.2 20 /

\plot 0 0 0 20 / \plot -0.2 0 -0.2 20 / \plot 0.2 0 0.2 20 /

\plot -20 0 -20 20 / \plot -20.2 0 -20.2 20 / \plot -19.8 0 -19.8 20 /

\plot -40 0 -40 20 / \plot -39.8 0 -39.8 20 / \plot -40.2 0 -40.2 20 /

\plot 40 0 20 20 /  
\plot 40 0 0 20 /  
\plot 40 0 -20 20 /  
\plot 40 0 -40 20 /  

\plot 20 0 0 20 /  
\plot 20 0 -60 20 /  

\plot 0 0 -20 20 /  
\plot 0 0 -40 20 /  

\plot -20 0 -40 20 /  
\plot -40 0 -60 20 /  

\endpicture 

\caption{$P_{(H,M)}(y_6):y_6x_5y_5x_4y_4x_3y_3$ and 
$n_2(y_6)=n_2(y_3)=2$}
\relabel{f5}
\end{figure}

\subsection{Externally B-active elements with respect to 
$M$}\relabel{sect5-2}

For any $y\in Y-Y'$, $y$ is the only vertex in the path
$P_{(H,M)}(y)$ belonging to $Y-Y'$.
We say $y$ is  
{\it externally B-active with respect to $M$} in $H$
\label{page-B-active}
if $w(y)>w(y_{\tau_{j_r}})$ holds for all $r=1,2,\cdots,t$,
where $\{y_{\tau_{j_r}}:r=1,2,\cdots,t\}
=Y'\cap V(P_{(H,M)}(y))$.
Let $A_{ex}(H,M)$ denote the set of those members 
in $Y-Y'$ 
which are externally B-active with respect to $M$ in $H$, and 
let $NA_{ex}(H,M)=(Y-Y')-A_{ex}(H,M)$.
Thus $NA_{ex}(H,M)$ is the set of those members 
in $Y-Y'$ 
which are not externally B-active with respect to $M$ in $H$.

Recall that the weight function 
$w$ of $G$ is a fixed injective mapping from $E$ to $\N_0$.
Introduced by Tutte \cite{tut}, 
for a given $T\in \sett(G)$,
an edge $y$ in $E(G)-E(T)$ is said to be 
{\it externally active with respect to $T$}
if $w(y)\ge w(y')$ holds for all edges $y'$ in the unique cycle of the subgraph $G[E(T)\cup \{y\}]$,
and an edge $y\in E(T)$ is said to be 
{\it internally active with respect to $T$} 
if $w(y)\ge w(y')$ holds for every edge $y'\in E(G)-E(T)$ 
with the property that 
$(E(T)-\{y\})\cup \{y'\}=E(T')$ holds 
for some $T'\in \sett(G)$. 
For the definition of these two concepts, 
the condition ``$w(y)\ge w(y')$"
can be replaced by ``$w(y)\le w(y')$",
as the condition is changed 
when $w(e)$ is replaced by $K-w(e)$ for each edge $e$ in $G$,
where $K$ is a number in $\N_0$ such that $K-w(e)\ge 0$ for all $e\in E$.
Tutte~\cite{tut} expressed the Tutte polynomial $T_G(x,y)$ 
as the summation of $x^{ia(T)}y^{ea(T)}$ 
over all spanning trees $T$ of $G$, where 
$ea(T)$ and $ia(T)$ are respectively the 
number of externally active edges and 
the number of internally active edges with respect to 
$T$.

In the following, we prove that 
the concept ``externally active with respect to $T$"
is extended to the one ``externally B-active 
with respect to $M$", where $M\in \UMX(H)$.

\begin{theo}\relabel{theo5-1}
Let $T\in \sett(G)$. For any
$y\in E(G)-E(T)$, 
$y$ is externally active respect to $T$ in $G$
\iff $y\in A_{ex}(H_{G,x_0},M_T)$.
\end{theo}

\proof Let $Y'=E(T)$ and let $H$ simply denote $H_{G,x_0}$
in the proof.
Thus $\sigma(H,Y')=1$,
and $\pi_i$'s and $\tau_i$'s are determined 
by Proposition~\ref{prop3-5}(i) and have the properties 
in Proposition~\ref{prop3-6}.

\def \setp {{\cal P}}

Write $x_{\pi_i}\preceq x_{\pi_j}$ 
if $x_{\pi_i}$ is a vertex on the path $P_{0,j}$
and $x_{\pi_i}\not\preceq  x_{\pi_j}$ otherwise. 
By Proposition~\ref{prop3-6} (iii), Claim 1 
follows directly.

\noindent {\bf Claim 1}: 
$x_{\pi_i}\preceq x_{\pi_j}$  implies that $i\le j$.

Thus  $x_{\pi_i}\preceq x_{\pi_j}$ \iff $i\le j$ and $P_{i,j}$ is part of $P_{0,j}$.
In the following, we first compare $i$ and $j$ 
in the case that 
$x_{\pi_i}\not \preceq x_{\pi_j}$ and 
$x_{\pi_j}\not \preceq x_{\pi_i}$.
Define $w_{max}(P_{i,j})$ as follows:
$$
w_{max}(P_{i,j})=
\left \{
\begin{array}{ll}
-1, &\mbox{if }E(P_{i,j})=\emptyset\\
\max \{w(e):  e\in E(P_{i,j})\},
\qquad &\mbox{otherwise}.
\end{array}
\right.
$$

\begin{figure}[htbp]
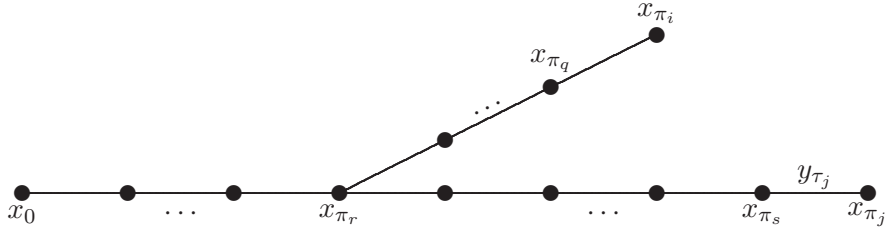

\centering

\beginpicture 
\setcoordinatesystem units <2pt,2pt> 
\setplotarea x from -110 to 80, y from 0 to 20 
\plotsymbolspacing=.2pt

\put {\cirs} at  -80 0
\put {\cirs} at  -60 0
\put {\cirs} at  -40 0
\put {\cirs} at  -20 0
\put {\cirs} at  0 0
\put {\cirs} at  20 0
\put {\cirs} at  40 0
\put {\cirs} at  60 0
\put {\cirs} at  80 0

\plot -80 0 80 0 /

\put {$x_0$} at  -80 -4
\put {$\cdots$} at  -50 -4
\put {$x_{\pi_r}$} at  -20 -4

\put {$\cdots$} at  30 -4
\put {$x_{\pi_s}$} at  60 -4
\put {$x_{\pi_j}$} at  80 -4
\put {$y_{\tau_j}$} at  70 3

\put {\cirs} at  40 30  \put {$x_{\pi_i}$} at  40 34  
\put {\cirs} at  20 20  \put {$x_{\pi_q}$} at  20 25  
\put {\cirs} at  0 10 
\plot -20 0 40 30 /
\put {$\cdot$} at  6 15
\put {$\cdot$} at  8 16
\put {$\cdot$} at  10 17

\endpicture 

\caption{
$x_{\pi_r}\preceq x_{\pi_i}$,
$x_{\pi_r}\preceq x_{\pi_j}$ but 
$E(P_{r,i})\cap E(P_{r,j})=\emptyset$.
}
\relabel{f8}
\end{figure}

\noindent {\bf Claim 2}: 
If $x_{\pi_r}\preceq x_{\pi_i}$,
$x_{\pi_r}\preceq x_{\pi_j}$ and
$E(P_{r,i})\cap E(P_{r,j})=\emptyset$,
then $w_{max}(P_{r,i})<w_{max}(P_{r,j})$ implies that 
$i<j$.

Assume that $w_{max}(P_{r,i})<w_{max}(P_{r,j})$.
We shall prove Claim 2 by induction on the 
the value of $\rho(i,j)=|E(P_{r,i})|+|E(P_{r,j})|$.
By Proposition~\ref{prop3-6} (iv) and the definition of $w_{max}(P_{i,j})$, 
Claim 2 holds when $|E(P_{r,i})|\le 1$ and $|E(P_{r,j})|\le 1$.

Assume that Claim 2 holds when $\rho(i,j)<K$,
where $K\ge 3$. 
Now consider the case that $\rho(i,j)=K$.

Let $k$ be the least possible integer such that 
$y_{\tau_k}$ is an edge on the path $P_{r,j}$
with $w(y_{\tau_k})>w_{max}(P_{r,i})$.
As $w_{max}(P_{r,i})<w_{max}(P_{r,j})$, such $k$ exists. 
By Claim 1, $r<k\le j$. 
If $k<j$, then $\rho(i,k)<K$ and by the inductive hypothesis, $w_{max}(P_{r,i})<w(y_{\tau_k})=w_{max}(P_{r,k})$
implies that $i<k$, and so $i<j$ holds.
Thus it suffices to consider the case that $k=j$,
i.e., $w_{max}(P_{r,i})<w(y_{\tau_j})$,
but $w_{max}(P_{r,i})>w(y_{\tau_t})$ 
for all edges $y_{\tau_t}$ on the path $P_{r,j}$ 
with $t\ne j$.

Let $s=b(y_{\tau_j})$ and $q=b(y_{\tau_i})$,
as shown in Figure~\ref{f8}, 
where $b(y_{\tau_j})$ is defined in Proposition~\ref{prop3-6}(iv)
(i.e., $b(y_{\tau_j})$ is the number $s$ 
such that $x_{\pi_s}$ is the end of $y_{\tau_j}$ in $G$ 
different from $x_{\pi_j}$). 
By Claim 1,  
$q<i$ and $s<j$. 
As $\rho(q,j)<K$, by the inductive hypothesis, 
$w(y_{\tau_j})>w_{max}(P_{r,i})\ge w_{max}(P_{r,q})$
implies that $j>q$.
As 
$w_{max}(P_{r,i})>w_{max}(P_{r,s})$, 
we have  $i>s$ by the inductive hypothesis.
Since $b(y_{\tau_j})=s<i$ and $b(y_{\tau_i})=q<j$,
the inequality $w(y_{\tau_j})>w_{max}(P_{r,i})\ge w(y_{\tau_i})$
implies that $j>i$ by 
Proposition~\ref{prop3-6} (iv).

Hence Claim 2 holds.

Now let $y$ be any edge in $E(G)-E(T)$. 
Assume that  $x_{\pi_i}$ 
and $x_{\pi_{j_1}}$ are the two ends of $y$, where $j_1>i$,
and the unique cycle $C$ in the graph obtained from $T$ by adding $y$ 
consists of edge $y$ and 
two edge-disjoint paths $P_{r,i}$ and $P_{r,j_1}$, 
where 
$x_{\pi_r}\preceq x_{\pi_i}$ and 
$x_{\pi_r}\preceq x_{\pi_{j_1}}$. 
Thus $r\le i<j_1$ with the possibility that $i=r$.

Let $x_{\pi_{j_1}}x_{\pi_{j_{2}}}\cdots x_{\pi_{j_t}}$  
be the longest possible subpath of $P_{r,j_1}$ 
between $x_{\pi_{j_1}}$ and $x_{\pi_{j_t}}$ 
such that $i<j_t$, as shown in Figure~\ref{f6}. 
By Claim 1, 
we have  
\begin{equation}\relabel{eq-4-2-2}
j_1>j_2>\cdots>j_t>i\ge b(y_{\tau_{j_t}}),
\end{equation}
where $i=b(y_{\tau_{j_t}})$ \iff 
$i=r$ and $b(y_{\tau_{j_t}})=r$.

\begin{figure}[htbp]
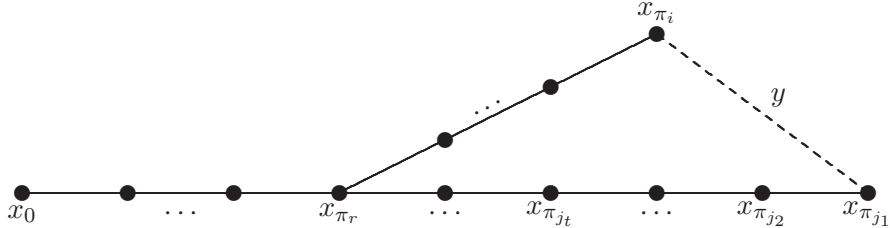

\centering

\beginpicture 
\setcoordinatesystem units <2pt,2pt> 
\setplotarea x from -110 to 80, y from 0 to 20 
\plotsymbolspacing=.2pt

\put {\cirs} at  -80 0
\put {\cirs} at  -60 0
\put {\cirs} at  -40 0
\put {\cirs} at  -20 0
\put {\cirs} at  0 0
\put {\cirs} at  20 0
\put {\cirs} at  40 0
\put {\cirs} at  60 0
\put {\cirs} at  80 0

\plot -80 0 80 0 /

\put {$x_0$} at  -80 -4
\put {$\cdots$} at  -50 -4
\put {$x_{\pi_r}$} at  -20 -4

\put {$\cdots$} at  0 -4
\put {$\cdots$} at  40 -4
\put {$x_{\pi_{j_t}}$} at  20 -4
\put {$x_{\pi_{j_2}}$} at  60 -4
\put {$x_{\pi_{j_1}}$} at  80 -4
\put {$y$} at  63 18

\put {\cirs} at  40 30  \put {$x_{\pi_i}$} at  40 34  
\put {\cirs} at  20 20  
\put {\cirs} at  0 10 
\plot -20 0 40 30 /
\put {$\cdot$} at  6 15
\put {$\cdot$} at  8 16
\put {$\cdot$} at  10 17

\linethickness=2pt 
\setdashes <0.1cm>
\plot 80 0 40 30 /

\endpicture 

\caption{$b(y_{\tau_{j_t}})\le i<j_t<\cdots<j_2<j_1$
}
\relabel{f6}
\end{figure}

As $j_1>i\ge b(y_{\tau_{j_t}})$,
by Claim 2, $w_{max}(P_{r,k})
\le w_{max}(P_{r,i})<w_{max}(P_{r,j_t})$, 
where $k=b(y_{\tau_{j_t}})$,
implying that 
$$
\max \{w_{max}(P_{r,i}), 
w_{max}(P_{r,j_1})\}
=\max\{w(y_{\tau_{j_s}}):s=1,2,\cdots,t\}.
$$
Thus the following claim holds.

\noindent {\bf Claim 3}: 
$y$ is externally active with respect to $T$ in $G$
\iff $w(y)>w(y_{\tau_{j_s}})$ holds for all 
$s=1,2,\cdots,t$.

On the other hand, 
by (\ref{eq-4-2-2})  and the fact that 
$n_1(y)=j_1$, $n_2(y_{\tau_{j_s}})=b(y_{\tau_{j_s}})=j_{s+1}$
for $s=1,2,\cdots,t-1$ and 
$n_2(y)=b(y)=i\ge k=b(y_{\tau_{j_t}})= n_2(y_{\tau_{j_t}})$,
the path $P_{M_T}(y)$ in $H_{G,x_0}$ with respect to $M_T$ 
is exactly the following one:
$$
P_{M_T}(y): yx_{\pi_{j_1}}y_{\tau_{j_1}} \cdots 
x_{\pi_{j_t}}y_{\tau_{j_t}}.
$$
Thus, by definition, the following claim also holds.

\noindent{\bf Claim 4}:
$y\in A_{ex}(H_{G,x_0},M_T)$ 
\iff $w(y)>w(y_{\tau_{j_s}})$ holds for all 
$s=1,2,\cdots,t$.

By Claims 3 and 4, the result holds.
\proofend

\subsection{Interpret B-parking function $f=\psi_H(M)$
\relabel{sect5-3}}

By the definition of the path $P_{(H,M)}(y)$, 
the following lemma follows.

\begin{lem}\relabel{lem5-1}
For any $y\in Y-Y'$, 
$y$ is adjacent to $x_{\pi_k}$ on the path $P_{(H,M)}(y)$
\iff 
$y\in N_H(x_{\pi_k})-\bigcup_{k<i\le n}N_H(x_{\pi_i})$.
\end{lem}

\begin{theo}\relabel{theo4-1} 
For any $y\in Y-Y'$ and $1\le k\le n$, 
$y\in D(H,Y',x_{\pi_k})$ 
\iff $y\in N_H(x_{\pi_k})-\bigcup_{k<i\le n}N_H(x_{\pi_i})$ and 
$y\in NA_{ex}(H,M)$. 
\end{theo}

\proof 
For $i=1,2,\cdots,n$, 
let $H_i$ be the subgraph of $H$ induced by 
$\sum_{i\le s\le n}N[x_{\pi_s}]$.
From Algorithm A, $\bigcup_{x\in X}D(H,Y',x)$ is a subset of $Y-Y'$ 
and can be partitioned into 
$n$ subsets $D'_1,D'_2,\cdots,D'_n$,
where $D'_i$ is the set of those vertices $y$ in $H_i$ 
having properties below:

\label{d-property}
(a) $y\in (Y-Y')-\bigcup_{1\le s<i}D'_i$;

(b) $y\in L(H_i)$;

(c) $w(y)<w(y_{\tau_i})$. 

Notice that $D'_i$ is the set of those members $y\in Y-Y'$ 
which are put into some set $D(x')$,
where $x'$ is the only neighbor of $y$ in $H_{i}$,
at Step A5 in Algorithm A
after $y_{\tau_{i-1}}$ is confirmed but before 
$y_{\tau_{i}}$ is confirmed. 

By Corollary~\ref{cor3-0}, 
if $y_{\tau_j}\in L(H_i)$, 
we have $j\ge i$ and $w(y_{\tau_j})\ge w(y_{\tau_i})$.
Thus the following claim holds:

\noindent {\bf Claim 1}: 
If $y_{\tau_j}\in L(H_i)$, 
then $w(y)<w(y_{\tau_j})$ holds for all $y\in D'_i$.

Now let $y$ be a member in $Y-Y'$.
Assume that $n_1(y)=j_1$ and 
the path $P_{(H,M)}(y)$ is 
$yx_{\pi_{j_1}}y_{\tau_{j_1}}\cdots 
x_{\pi_{j_t}}y_{\tau_{j_t}}$.
By the definition of  $P_{(H,M)}(y)$,
$j_{s+1}=n_2(y_{\tau_{j_{s}}})$
for $s=1,2,\cdots,t-1$ 
and $j_1>j_2>\cdots >j_t>n_2(y)\ge n_2(y_{\tau_{j_{t}}})$.

($\Rightarrow$)
Assume that $y\in D(H,Y',x_{\pi_k})$.
By Proposition~\ref{prop3-1} (ii), 
$y\in N_H(x_{\pi_k})$.
Let $i$ be the minimum integer with $0<i\le k$
such that $y\in L(H_{i})$
and $w(y)<w(y_{\tau_{i}})$.
Such $i$ exists by Proposition~\ref{prop3-1} (ii).
Thus $y\in D_{i}'$.
Clearly, $k=n_1(y)=j_1$ and so 
$x_{\pi_k}$ (i.e., $x_{\pi_{j_1}}$) is the vertex on the 
path $P_{(H,M)}(y)$ adjacent to $y$.
It remains to show that $y\in  NA_{ex}(H,M)$.

As  $y\in L(H_i)$, 
we have $q<i\le j_1$, where $q=n_2(y)$.
Note that $j_1>j_2>\cdots>j_t>n_2(y)=q\ge n_2(y_{\tau_{j_{t}}})$.
Thus $j_{s+1}<i\le j_s$ holds for some $s$ with $1\le s\le t$,
where assume that $j_{t+1}=n_2(y)=q$ when $s=t$.
Then $y_{\tau_{j_s}}\in L(H_i)$.
By Claim 1, $w(y)<w(y_{\tau_{j_s}})$. 
By definition, $y\in  NA_{ex}(H,M)$.

Hence the necessity holds.

($\Leftarrow$) 
Now assume that 
$y\in  NA_{ex}(H,M)$.
Assume that $j_1=n_1(y)$. 
We will show that $y\in D(H,Y',x_{\pi_{j_1}})$.

On the contrary, suppose that 
$y\notin D(H,Y',x_{\pi_{j_1}})$.
By Proposition~\ref{prop3-1} (ii), 
$y\notin D(H,Y',x_{\pi_{s}})$
for all $s=1,2,\cdots,n$, implying that  
$y\notin D'_i$ for all $i=1,2,\cdots,n$.

As $j_1=n_1(y)$ and $q=n_2(y)$,  
$y\in L(H_i)$ for all $i$ with $q<i\le {j_1}$.
For each $i$ with $q<i\le {j_1}$,
as $y\notin D'_i$,
we have $w(y)>w(y_{\tau_i})$ by property (c).
Particularly,  as $q<j_t<\cdots<j_1$, 
$w(y)>w(y_{\tau_{j_s}})$ holds 
for all $s=1,2,\cdots,t$,
implying that $y$ is externally B-active
with respect to $M$ in $H$.
Thus $y\not\in  NA_{ex}(H,M)$,
a contradiction.

Hence the sufficiency holds.
\proofend

By Theorem~\ref{theo4-1} and the definition of $\psi_H$, 
we have the following 
corollaries. 

\begin{cor}\relabel{cor4-1} 
Let $M\in \UMX(H)$. 
If $f=\psi_H(M)$, then,  
$f(x_{\pi_i})$ is the size of the set 
$\left (N_H(x_{\pi_i})\cap NA_{ex}(H,M)\right )
-\bigcup_{i<s\le n}N_H(x_{\pi_s})$
for all $i=1,2,\cdots,n$.
\end{cor}

\begin{cor}\relabel{cor4-3}
Let $M\in \UMX(H)$. If $f=\psi_H(M)$, then
$$
\sum_{x\in X}f(x)=|NA_{ex}(H,M)|\le |Y|-|X|.
$$ 
\end{cor}

Now we apply Theorem~\ref{theo5-1}
to find another interpretation for 
G-parking functions of $G$.

Let $T\in \sett(G)$.
Write $x_{\pi_i}\ll_T x_{\pi_j}$ 
for all $i,j$ with $0\le i<j\le n$.
For any two vertices $x'$ and $x$ in $G$, let 
$P_T(x',x)$ denote the unique path in $T$ between $x'$ and $x$.

\begin{pro}\relabel{prop5-2}
For any two different vertices $x'$ and $x$ in $G$, 
the following statements are equivalent:
\begin{enumerate}
\item $x'\ll_T x$ ;
\item 
$w_{max}(P_T(x'',x'))<w_{max}(P_T(x'',x))$,
where $x''$ is the vertex in both paths $P_T(x_0,x')$ and $P_T(x_0,x)$ 
with $E(P_T(x'',x'))\cap E(P_T(x'',x))=\emptyset$;
\item if $y$ is an edge in $E(G)-E(T)$ joining $x$ and $x'$, 
then $x$ is the vertex $x_{\pi_j}$ with $j=n_1(y)$,
where $\pi_s=\pi_s(H_{G,x_0},E(T))$ for $s\in \{1,2,\cdots,n\}$; 
\item 
if $y$ is an edge in $E(G)-E(T)$ joining $x$ and $x'$, then 
$x$ is the vertex in the path $P_{(H,M)}(y)$ adjacent to $y$,
where $Y'=E(T)$.
\end{enumerate}
\end{pro}

\proof Claims 1 and 2 in the proof of Theorem~\ref{theo5-1} imply that 
(i) $\Leftrightarrow $ (ii), while the definition of the path 
$P_{(H,M)}(y)$ implies that (iii) $\Leftrightarrow $ (iv).
Finally, by the definition of the ordering $\ll_T$ and 
the definition of $n_1(y)$, (i) $\Leftrightarrow$ (iii) follows. 
\proofend

Recall that the mapping 
$\phi_G: \sett(G)\rightarrow \GP(G,x_0)$ is 
defined by $\phi_G(T)=\psi_{H_{G,x_0}}(M_T)$,
where $M_T=\{x_{\pi_i}y_{\tau_i}:i=1,2,\cdots,n\}$
by Corollary~\ref{cor3-2}.
By Corollary~\ref{cor4-1} and Proposition~\ref{prop5-2}, 
we get the following interpretation for $\phi_G$ 
which is different from the one in Corollary~\ref{cor3-3}.

\begin{cor}\relabel{cor4-4} 
Let $T\in \sett(G)$.
If $f=\phi_G(T)$, then, for any $x\in V-\{x_0\}$, 
$f(x)$ is the number of those edges $y\in E(G)-E(T)$ 
such that 
$y$ is not externally active with respect to $T$ in $G$
and
$y$ is incident with $x$ and $x'$, where $x'\ll_T  x$.
\end{cor}

By Corollaries~\ref{cor4-3} and~\ref{cor4-4}, 
we have the following conclusion.

\begin{cor}\relabel{cor4-5}
Let $T\in \sett(G)$.
If $f=\phi_G(T)$,  
then  
$$
ea(T)+\sum_{x\in X}f(x)=|E(G)|-|V(G)|+1,
$$
where $ea(T)$ is the number of 
externally active edges with respect to $T$ in $G$.
\end{cor}

\subsection{A generating function $\Omega(H;x,y,z)$ 
\relabel{gen-func}}

Let $M\in \UMX(H)$.
For any $x_{\pi_q}\in X$, let $R(x_{\pi_{q}})$ 
denote the following unique path:
$$
x_{\pi_{j_1}} y_{\tau_{j_1}}x_{\pi_{j_2}}y_{\tau_{j_2}}
\cdots x_{\pi_{j_s}}y_{\tau_{j_s}},
$$ 
where $j_1=q$, 
$j_{i+1}=n_2(y_{\tau_{j_i}})$ for $i=1,2,\cdots,s-1$
and $y_{\tau_{j_s}}\in L(H)$.
For any $y'\in Y-V(M)$ and $r\ge 1$,
if $t_r=n_r(y')\ge 1$, 
let $Q_r(y')$ be the path in $H$ formed by 
combining edge $y'x_{\pi_{t_r}}$ and path $R(x_{\pi_{t_r}})$.
In the case that $y'\in L(H)$ (i.e., $n_2(y')<1$),
assume that $Q_2(y')$ consists of vertex $y'$ only.
Let $k=0$ if $V(Q_1(y'))\cap V(Q_2(y'))
\cap X=\emptyset$, and let $k$ be the largest integer 
with $x_{\pi_k}\in V(Q_1(y'))\cap V(Q_2(y'))
\cap X$ otherwise. 
Let $C_H(y')$ be the set 
$
\{y_{\tau_u}\in V(Q_1(y'))\cup V(Q_2(y')): u>k, y_{\tau_u}\ne y'\}.
$

For example, if $H$ is the graph in Figure~\ref{f5}
and $M=\{x_iy_i:i=1,2,\cdots,5\}$,
then $Q_1(y_6)$ is the path 
$y_6x_5y_5x_4y_4x_3y_3x_2y_2x_1y_1$ and 
$Q_2(y_6)$ is the path 
$y_6x_2y_2x_1y_1$.
Thus $k=2$ and $C_H(y_6)=\{y_5,y_4,y_3\}$.
For the bipartite graph $H_{G,x_0}$
and $M=M_T$, where $T\in \sett(G)$,
$C_{H_{G,x_0}}(y')$ corresponds to the set of edges 
$y\ne y'$
in the unique cycle of $G[E(T)\cup \{y'\}]$,
where $y'\in E(G)-E(T)$.

For any $y_i\in V(M)\cap Y$,
$y_i$ is said to be {\it internally B-active 
with respect to $M$} if $w(y_i)>w(y')$ holds 
for each $y'\in Y-V(M)$ with $y_i\in C(y')$.
Let $A_{in}(H,M)$ be the set of 
internally B-active members 
with respect to $M$ in $H$.

Define a function $\Omega(H;x,y,z)$ with three
variable $x,y,z$ 
as follows:
\begin{equation}
\Omega(H;x,y,z)=
\sum_{S\subseteq X}z^{|X|-|S|}
\sum_{M\in \UMS(H)}x^{ia_S(M)}y^{ea_S(M)},
\end{equation}
where $ia_S(M)=|A_{in}(H[N[S]],M)|$ and 
$ea_S(M)=|A_{ex}(H[N[S]],M)|$.

If $\Omega(H;1,1,z)=\sum_{i\ge 0} c_iz^i$, then 
$c_i$ is the number of members $M\in \UM(H)$
with $|M|=|X|-i$.
In particular, $c_0=|\UMX(H)|$.

If $\Omega(H;x,y,0)=\sum_{i,j\ge 0} u_{i,j}x^iy^j$, then 
$u_{i,j}$ is the number of members $M\in \UMX(H)$
with $|A_{in}(H,M)|=i$  and $|A_{ex}(H,M)|=j$
(i.e., $|NA_{ex}(H,M)|=|Y|-|X|-j$).

If $\Omega(H;1,y,0)=\sum_{j\ge 0} d_jy^j$, then 
$d_j$ is the number of members $M\in \UMX(H)$
with $|A_{ex}(H,M)|=j$, i.e., 
$|NA_{ex}(H,M)|=|Y|-|X|-j$. 
By Corollary~\ref{cor4-3}, $d_j$ is the number 
of members $f$ in $\BPX(H)$ with 
$\sum_{x\in X}f(x)=|Y|-|X|-j$.
By Corollary~\ref{cor4-5}, if $H=H_{G,x_0}$, then 
$d_j$ is the number of members $f\in \GP(G,x_0)$ with 
$\sum_{x\in X}f(x)=|E(G)|-|V(G)|+1-j$.

For any $T\in \sett(G)$, an edge $e\in E(T)$  
is internally active with respect to $T$ in $G$ 
\iff $e\in A_{in}(H,M_T)$.
Thus, $\Omega(H_{G,x_0};x,y,0)$ is the Tutte 
polynomial $T_G(x,y)$.

\vspace{0.4 cm}

\noindent {\bf Acknowledgment}: 
The author wishes to thank the referees 
for their very helpful suggestions.

\end{document}